\documentclass{article}

\usepackage[english]{babel}

\usepackage[letterpaper,top=2cm,bottom=2cm,left=3cm,right=3cm,marginparwidth=1.75cm]{geometry}

\usepackage{authblk}
\usepackage{amsmath}
\usepackage{amssymb}
\usepackage{amsthm}
\usepackage{mathptmx}
\usepackage{caption}
\usepackage{graphicx}
\usepackage{subfigure}
\usepackage[colorlinks=true, allcolors=blue]{hyperref}
\newtheorem{definition}{Definition}[section]
\newtheorem{theorem}{Theorem}[section]

\newtheorem{lemma}{Lemma}[section]
\newtheorem{remark}{Remark}[section]
\newtheorem{corollary}{Corollary}[section]

\newtheorem{assumption}{Assumption}
\numberwithin{equation}{section} 
\makeatletter

\newcommand{\Rmnum}[1]{\expandafter\@slowromancap\romannumeral #1@}

\providecommand{\keywords}[1]
{
  \small	
  \textbf{Keywords:} #1
} 

\author[1]{Xu Cheng}
\author[2,3]{Jiaqi Liu}
\author[4,5]{Zaijiu Shang\thanks{Corresponding author, Email:zaijiu@simis.cn }}

\affil[1]{Naval University of Engineering, Wuhan 430033, China}
\affil[2]{Academy of Mathematics and Systems Science, Chinese Academy of Sciences, Beijing 100190, China}
\affil[3]{University of Chinese Academy of Sciences, Beijing 100049, China}
\affil[4]{Center for Mathematics and Interdisciplinary Sciences at Fudan University, Shanghai 200433, China}
\affil[5]{Shanghai Institute for Mathematics and Interdisciplinary Sciences, Shanghai 200438, China}
\makeatother
\title{A class of generalized Nesterov's accelerated gradient method from dynamical perspective\footnotemark[2]}
\graphicspath{{figures/}}

\begin{document}
\maketitle
\renewcommand{\thefootnote}{\fnsymbol{footnote}}
\footnotetext[2]{This work was supported by the National Natural Science Foundation of China (Grant No.12071471).}
\begin{abstract}
We propose a class of \textit{Euler-Lagrange} equations indexed by a pair of parameters ($\alpha,r$)  that generalizes Nesterov's accelerated gradient methods for convex ($\alpha=1$) and strongly convex ($\alpha=0$) functions from a continuous-time perspective. This class of equations also serves as an interpolation between the two Nesterov's schemes. The corresponding \textit{Hamiltonian} systems can be integrated via the symplectic Euler scheme with a fixed step-size. Furthermore, we can obtain the convergence rates for these equations ($0<\alpha<1$)  that outperform Nesterov's when time is sufficiently large for $\mu$-strongly convex functions,  without requiring a priori knowledge of $\mu$. We demonstrate this by constructing a class of Lyapunov functions that also provide a unified framework for Nesterov's schemes for convex and strongly convex functions.\\
\\
\keywords{convex optimization, accelerated gradient method, Euler-Lagrange equation, non-autonomous Hamiltonian system, splitting method, Lyapunov analysis}
\end{abstract}

\section{Introduction}
In this paper, we mainly focus on the unconstrained optimization problem:
\begin{equation}\label{prob}
    \mathop{\min}_{x\in\mathbb{R}^d} f(x)
\end{equation}
where $f:\mathbb{R}^d\to \mathbb{R}$ is a twice continuously differentiable \textit{convex} function and $\nabla f$ is Lipschitz continuous with Lipschitz constant $L$. We assume that $f$ has a unique minimizer, let $x^*=\mathop{\arg \min}_{x\in\mathbb{R}^d} f(x)$ be the minimizer and $f^*=f(x^*)=\mathop{\min}_{x\in\mathbb{R}^d} f(x)$. And we consider the first-order methods to solve 
the problem (\ref{prob}), which only involve
the first-order derivative of $f$. These methods are popular when the second and higher-order derivatives are difficult to compute, such as in machine learning and other large-scale problems.\par
The simplest first-order method is the gradient descent method (GD):
\begin{equation}\label{alg:GD}
    x_{n+1}=x_n-h\nabla f(x_n)
\end{equation}
where the step-size $h\leq \frac{1}{L}$, and its convergence rate is:
\begin{equation}
    f(x_n)-f^*\leq O(\frac{1}{n}).
\end{equation}
A milestone in first-order methods came from Nesterov, who constructed an \textit{accelerated} gradient method (NAG) in \cite{nesterov1983method} by adding a momentum term:
\begin{equation}\label{alg:NAG}
\begin{split}
        x_{n+1}&=y_n-h\nabla f(y_n)\\
        y_{n+1}&=x_{n+1}+\frac{n}{n+3}(x_{n+1}-x_n)
    \end{split}
\end{equation}
where the step-size $h\leq \frac{1}{L}$, and proved that it exhibits a faster convergence rate than GD:
\begin{equation}\label{convergence:NAG}
    f(x_n)-f^*\leq O(\frac{1}{n^2}).
\end{equation}
For \textit{$\mu$-strongly convex} functions, Nesterov \cite{nesterov2004introductory} introduced another scheme (NAG-SC):
\begin{equation}\label{alg:NAG-sc}
    \begin{split}
        x_{n+1}&=y_n-h\nabla f(y_n)\\
        y_{n+1}&=x_{n+1}+\frac{1-\sqrt{\mu h}}{1+\sqrt{\mu h}}(x_{n+1}-x_n)
    \end{split}
\end{equation}
and gave the convergence rate:
\begin{equation}
    f(x_n)-f^*\leq O(e^{-\sqrt{\kappa}n}),
\end{equation}
while the convergence rate of GD for \textit{$\mu$-strongly convex} functions is 
\begin{equation}
     f(x_n)-f^*\leq O(e^{-\kappa n})
\end{equation}
where $\kappa=\frac{\mu}{L}\leq 1$. Therefore, NAG-SC improves GD by a factor of $\sqrt{\kappa}$. It is worth mentioning that NAG-SC requires us to know $\mu$ a priori, which is hard in practice.\par
However, there remains a lack of understanding of the acceleration phenomenon and a framework to generate accelerated method for many years.
A significant progress in this regard was the derivation by Su et al. \cite{NAGODE} of NAG's limiting ordinary differential equation (\ref{alg:NAG}):
\begin{equation}\label{eqn:NAG}
    \Ddot{x}+\frac{3}{t}\Dot{x}+\nabla f(x)=0
\end{equation}
where $\cdot=\frac{d}{dt}$. Moreover, if $x(t)$ is the solution of (\ref{eqn:NAG}) with initial values $x(0)=x_0,\ \Dot{x}(0)=0$, they verified the convergence rate:
\begin{equation}
    f(x(t))-f^*\leq O(\frac{1}{t^2})
\end{equation}
using Lyapunov analysis, which is consistent with the discrete case (\ref{convergence:NAG}) and gave a class of generalized accelerated gradient descent methods when the constant "3" in (\ref{eqn:NAG}) is replaced by $r\geq 3$.\par
Since then, there have been many studies on optimization from a continuous-time perspective. Some of them derived continuous analog of a existing algorithm and exploited it to interpret some intriguing characteristics of the algorithm (e.g., \cite{shi2022understanding}, \cite{luo2022differential},  \cite{chen2022gradient}, \cite{li2024linear}). For example, Shi et al. \cite{shi2022understanding} derived a high-resolution differential equation of NAG (\ref{alg:NAG}):
\begin{equation}\label{eqn:NAG_highresolution}
    \Ddot{x}+\frac{3}{t}\Dot{x}+\sqrt{s}\nabla^2f(x)\Dot{x}+(1+\frac{3\sqrt{s}}{2t})\nabla f(x)=0
\end{equation}
which is a more accurate continuous analog of NAG.
Based on it, Shi et al. \cite{li2024linear} established the linear convergence of NAG (\ref{alg:NAG}) for \textit{$\mu$-strongly convex} functions without the knowledge of $\mu$:
\begin{equation}\label{convergence:NAG_linear}
    f(x(t))-f^*\leq \frac{\mathcal{E}(T)}{2t(t+\sqrt{s})}e^{-\mu\sqrt{s}(t-T)}
\end{equation}
when $t\geq T=\frac{4}{\mu\sqrt{s}}$, $x(t)$ is a solution of equation (\ref{eqn:NAG_highresolution}) and $\mathcal{E}(t)$ is a Lyapunov function of the equation.\par
Others attempted to generalize NAG from a variational perspective (e.g., \cite{1603}, \cite{betancourt2018symplectic}, \cite{wilson2021lyapunov}, \cite{diakonikolas2021generalized}, \cite{kim2023unifying}). Notably, Wibisono et al. \cite{1603} introduced the \textit{Bregman Lagrangian}
\begin{equation}\label{Lagrangian:Bregman}
    \mathcal{L}(x,v,t)=e^{\alpha(t)+\gamma(t)}[D_h(x+e^{-\alpha(t)}v,x)-e^{\beta(t)}f(x)]
\end{equation}
where $D_h$ is the \textit{Bregman Divergence} with respect to a convex and essentially smooth function $h:\mathbb{R}^d\to \mathbb{R}$ defined as follows:
\begin{equation*}
    D_h(y,x)=h(y)-h(x)-<\nabla h(x),y-x>
\end{equation*}
and $\alpha, \beta,\gamma:\mathbb{R_+}\to\mathbb{R}$ are continuously differentiable functions with respect to time. We say it is in the \textit{Euclidean setting} if $h(x)=\frac{1}{2}||x||^2$. They proved that for \textit{convex} function $f$, if $\alpha, \beta,\gamma$ satisfy the \textit{ideal scaling conditions}
\begin{equation}\label{conditions:ideal}
    \Dot{\beta}\leq e^{\alpha},\
    \Dot{\gamma}=e^{\alpha},
\end{equation}
then the solutions of the \textit{Euler-Lagrange} equation
\begin{equation}\label{eqn:E-L_Bregman}
    \frac{d}{dt}\frac{\partial \mathcal{L}}{\partial v}(x,\Dot{x},t)=\frac{\partial \mathcal{L}}{\partial x}(x,\Dot{x},t)
\end{equation}
satisfy
\begin{equation}\label{convergence:Bregman_general}
    f(x(t))-f^*\leq O(e^{-\beta(t)}).
\end{equation}
In particular, they defined the polynomial subfamily of the \textit{Bregman Lagrangians} with the following choices of $\alpha, \beta,\gamma$ indexed by a parameter $p$:
\begin{equation}\label{conditions:polynomial}
\alpha(t)=\ln p-\ln t,\ \beta(t)=p\ln t+\ln C,\ \gamma(t)=p\ln t
\end{equation}
that attains a polynomial convergence rate of degree $p$.
Indeed, with the choices (\ref{conditions:polynomial}), $\alpha, \beta,\gamma$ satisfy the \textit{ideal scaling conditions} (\ref{conditions:ideal}), and thus the objective function converges to its minimum at a rate of  $O(t^{-p})$ due to  (\ref{convergence:Bregman_general}). In particular, if we take $p=2, C=
\frac{1}{4}$ in the \textit{Euclidean setting}, then the \textit{Euler-Lagrange} equation (\ref{eqn:E-L_Bregman}) is exactly the NAG's ODE (\ref{eqn:NAG}). Therefore, the \textit{Bregman} dynamical system can be viewed as a generalization of (\ref{eqn:NAG}). 
\par
However, the convergence rate $O(t^{-p})$ can not be obtained easily via a single-loop algorithm like NAG (\ref{alg:NAG}). In order to illustrate that, we consider the following \textit{time-dependent Lagrangian}
\begin{equation}\label{Lagrangian:L_ab}
    \mathcal{L}(x,v,t)=\frac{1}{2}a(t)||v||^2-b(t)f(x)
\end{equation}
with positive time-dependent coefficients $a(t)$ and $b(t)$ and $t\geq0$. The \textit{Bregman Lagrangian} in the \textit{Euclidean setting} is a special case of (\ref{Lagrangian:L_ab}) when we take $a(t)=e^{\gamma(t)-\alpha(t)}$ and $b(t)=e^{\alpha(t)+\beta(t)+\gamma(t)}$.
Then the corresponding \textit{Euler-Lagrange} equation is
\begin{equation}
    \Ddot{x}+\frac{\Dot{a}}{a}\Dot{x}+\frac{b}{a}\nabla f(x)=0.
\end{equation}
We can convert it into its equivalent \textit{Hamiltonian} system via \textit{Legendre} transformation $y=\partial_v \mathcal{L}=a(t)v$:
\begin{equation}
\begin{cases}
    \Dot{x}&=\frac{1}{a(t)}y\\
    \Dot{y}&=-b(t)\nabla f(x).\\
\end{cases}
\end{equation}
Then we can apply the symplectic Euler method
\begin{equation}\label{alg:SE_ab}
    \begin{split}
        y_{n+1}&=y_n-h_nb(t_n)\nabla f(x_n)\\
        x_{n+1}&=x_n+h_n\frac{1}{a(t_n)}y_{n+1}\\
        t_{n+1}&=t_n+h_n
    \end{split}
\end{equation}
to integrate the system. And we compute the linear stability region of the scheme:
\begin{theorem}\label{thm1}
    The symplectic Euler method (\ref{alg:SE_ab}) is (linearly) stable if the step-size sequence $\{ h_n\}_{n=0}^\infty$ satisfies
    \begin{equation}\label{conditions:stepsize}
        h_n^2\frac{b(t_n)}{a(t_n)}L<4,\ n=0,1,2...
    \end{equation}
    where $t_n=t_0+\sum_{k=0}^{n-1}h_k$, $L$ is the Lipschitz constant of $\nabla f$.
\end{theorem}

\begin{corollary}\label{cor1}
    For the polynomial subfamily of the \textit{Bregman Lagrangians} in the \textit{Euclidean setting}, the algorithm (\ref{alg:SE_ab}) is stable if the step-size sequence $\{ h_n\}_{n=0}^\infty$ satisfies
    \begin{equation}\label{conditions:stepsize_poly}
        h_n<\sqrt{\frac{4}{CLp^2(t_0+\sum_{k=0}^{n-1}h_k)^{p-2}}},\ n=0,1,2...
    \end{equation}
    In particular, when $p=2$ (i.e.NAG case), we can have a fixed step-size
    \begin{equation}
        h<\sqrt{\frac{1}{CL}}.
    \end{equation}
\end{corollary}
\begin{remark}
    (i) If $p>2$, then $h_n\to 0$ as $n\to \infty$.\\
    (ii) We consider the case where $p=4$. Let $t_0=1,\ C=\frac{1}{Lp^2}$, then we can choose $h_n=\frac{1}{t_n}$ satisfying the stability condition (\ref{conditions:stepsize_poly}). Therefore, $t_{n+1}=t_n+h_n=t_n+\frac{1}{t_n},\ t_0=1$. Then we can approximate $t_n\sim\sqrt{2n}$ as $n\to \infty$. From (\ref{convergence:Bregman_general}) and (\ref{conditions:polynomial}), we know that $f(x_n)-f^*\approx O(t_n^{-4})\approx O((\sqrt{2n})^{-4})\approx O(n^{-2})$. It seems unlike to attain  O($n^{-4}$) in this case.\end{remark} 
    Motivated by this, we attempt to construct dynamical systems which have faster convergence rate than NAG and can be numerically solved by a single-loop algorithm with a fixed step-size. We consider the following \textit{Lagrangian}
\begin{equation}\label{Lagrangian:L_phi}
    \mathcal{L}(x,v,t)=e^{\xi(t)}[\frac{1}{2}||v||^2-f(x)]
\end{equation}
which is a special case of (\ref{Lagrangian:L_ab}) with $a(t)=b(t)=e^{\xi(t)}$. Therefore, we know that its \textit{Euler-Lagrange} equation can be integrated by a fixed step-size symplectic Euler method from Theorem\ref{thm1}. The \textit{Euler-Lagrange} equation of (\ref{Lagrangian:L_phi}) is
\begin{equation}
    \Ddot{x}+\Dot{\xi}\Dot{x}+\nabla f(x)=0.
\end{equation}
We focus on the case $\Dot{\xi}(t)=\frac{r}{t^\alpha}$, where $r>0$ and $0\leq \alpha \leq 1$, then the \textit{Euler-Lagrange} equation becomes
\begin{equation}\label{eqn:E-L_ar}
    \Ddot{x}+\frac{r}{t^\alpha}\Dot{x}+\nabla f(x)=0.
\end{equation}
We call it ($\alpha,r$)-damped system. When $\alpha=1,\ r=3$, it is the continuous limit of NAG (\ref{alg:NAG}). If $\alpha=0,\ r=2\sqrt{\mu}$, it is the continuous limit of NAG-SC (\ref{alg:NAG-sc}) and \cite{wilson2021lyapunov} introduced another \textit{Bregman Lagrangian} that generalizes it but can not recover the above \textit{Bregman Lagrangian} (\ref{Lagrangian:Bregman}). \cite{kim2023unifying} proposed a unified  \textit{Bregman Lagrangian} to address the inconsistency between the two types of  \textit{Bregman Lagrangians}. When $\alpha=1$, \cite{NAGODE} gave some results when $1<r<3$ (low friction) and $r>3$ (high friction). To the best of our knowledge, there is few research on the case where $0<\alpha<1$. In this paper, we provide a unified framework as an interpolation between the convex and strongly convex cases.
\begin{remark}
    We consider the linearization of equation (\ref{eqn:E-L_ar}) i.e. $f(x)=\sum_{i=1}^d\frac{1}{2}\lambda_ix_i^2$, then\\
    (i) for $\alpha=0$, equation (\ref{eqn:E-L_ar}) is 
    \begin{equation}
        \Ddot{x}+r\Dot{x}+\lambda x=0
    \end{equation}
    which is a second-order ODE with constant coefficients.\\
    (ii) for $\alpha=1$, equation (\ref{eqn:E-L_ar}) is 
    \begin{equation}
        \Ddot{x}+\frac{r}{t}\Dot{x}+\lambda x=0
    \end{equation}
     which can be transformed to the Bessel equation of order $\frac{r-1}{2}$ via $y(\tau)={\tau}^{\frac{r-1}{2}}x(t(\tau)),\ \tau=\sqrt{\lambda}t$.\\
     (iii) for $0<\alpha<1$, equation (\ref{eqn:E-L_ar}) is not equivalent to a well-studied ODE and $t=0$ is not a regular singular point of it since $t^\alpha$ is not analytic. Therefore, we cannot solve it by power series method.
\end{remark}
Equation (\ref{eqn:E-L_ar}) is a dissipative system with the dissipative term $\frac{r}{t^\alpha}\Dot{x}$. The dissipative coefficient $\frac{r}{t^\alpha}\to 0$ as $t\to \infty$, and it converges faster as $\alpha$ increases. Intuitively, the system may dissipate faster as $\alpha$ gets smaller when $t$ is large considering that the dissipative term disappears more slowly. Therefore, we hope that $f$ decreases faster along the solution when $\alpha<1$ than $\alpha=1$ (NAG) as $t$ is large enough. Our main results are as follows:
\begin{theorem}\label{thm:Lyapunov}
    Assume $f(x)$ is a $\mu$-strongly convex function, then
    \begin{equation}
        L(t)=\frac{1}{2}||A(t)\Dot{x}(t)+B(t)(x(t)-x^*)||^2+\frac{1}{2}C(t)||\Dot{x}(t)||^2+D(t)(f(x(t))-f^*)
    \end{equation}
    is a Lyapunov function of the equation (\ref{eqn:E-L_ar}) when $t\geq T=T(\mu,\alpha,r,
    \beta)$, where $T(\mu,\alpha,r,
    \beta)$ is either the largest root of the polynomial $P(t;\mu,\alpha,r,\beta)=\mu t^2-(1-\beta)\beta r^2t^{2(1-\alpha)}+(2-3\beta)r\alpha t^{1-\alpha}+\alpha(\alpha+1)$ or zero if $P(t;\mu,\alpha,r,\beta)$ has no root and
    \begin{align}
        D(t)&=e^{\beta\xi(t)}\label{expression:Dt}\\
        B(t)&=\sqrt{\beta r[(1-\beta)rt^{1-\alpha}+\alpha]t^{-\alpha-1}e^{\beta\xi(t)}}\label{expression:Bt}\\
        A(t)&=\sqrt{\frac{\beta rt^{1-\alpha}}{(1-\beta)rt^{1-\alpha}+\alpha}e^{\beta\xi(t)}}\label{expression:At}\\
        C(t)&=[1-\frac{\beta rt^{1-\alpha}}{(1-\beta)rt^{1-\alpha}+\alpha}]e^{\beta\xi(t)}\label{expression:Ct}
    \end{align}
    with 
    \begin{equation}\label{conditions:beta_r}
        \xi(t)=\begin{cases}
            \frac{r}{1-\alpha}t^{1-\alpha},\ &0\leq \alpha<1\\
            r\ln t,\ &\alpha=1
        \end{cases}
        ,\ 
            0\leq\beta\leq \begin{cases}
            \frac{1}{2}\quad,\ &0\leq \alpha<1\\
            \min\{\frac{2}{3},\frac{1+r}{2r}\},\ &\alpha=1
        \end{cases},\ r\leq\sqrt{\frac{\mu}{\beta(1-\beta)}}\ when\ \alpha=0.
    \end{equation}
\end{theorem}
\begin{theorem}\label{thm:convergence_ar}
    Assume $f(x)$ is a $\mu$-strongly convex function, $x(t)$ is a solution of the equation (\ref{eqn:E-L_ar}) with initial conditions $x(t_0)=x_0,\ \Dot{x}(t_0)=0$, then\\
    (i) $\alpha=0$, $f(x(t))-f^*\leq[\frac{\mu}{2}||x_0-x^*||^2+f(x_0)-f^*]e^{-\sqrt{\mu}t}$ when $t\geq t_0=0$ and $r=2\sqrt{\mu}$.\\
    (ii) $0<\alpha<1$, $f(x(t))-f^*\leq L(T)e^{-\frac{r}{2(1-\alpha)}t^{1-\alpha}}$, when $t\geq T$, $T$ and $L(t)$ are defined in Theorem \ref{thm:Lyapunov} with $\beta=\frac{1}{2}$.\\
    (iii) $\alpha=1$, $f(x(t))-f^*\leq L(T)t^{-\frac{1+r}{2}}$, when $r>3$, $t\geq T=\sqrt{\frac{(r-3)(r+1)}{4\mu}}$ and L(t) is defined in Theorem \ref{thm:Lyapunov} with $\beta=\frac{1+r}{2r}$, or $f(x(t))-f^*\leq [\frac{r}{3}(\frac{r}{3}+1)t_0^{\frac{2}{3}r-2}||x_0-x^*||^2+t_0^{\frac{2}{3}r}(f(x_0)-f^*)]t^{-\frac{2}{3}r}$, when $t\geq t_0$ and $r\leq 3$.
    \end{theorem}
\begin{remark}
    In \cite{li2024linear}, they proved linear convergence (\ref{convergence:NAG_linear}) of NAG for $\mu$-strongly convex functions. In contrast, we obtain a faster convergence rate than NAG when $\mu$ is small enough and $t$ is large enough in Theorem \ref{thm:convergence_ar} (ii).
\end{remark}
The remainder of the paper is organized as follows. In Section \ref{sect:2}, we introduce symplectic integrators of non-autonomous \textit{Hamiltonian} systems via splitting methods. The proof of Theorem \ref{thm1} is shown in this section. In Section \ref{sect:3}, we construct Lyapunov functions and prove Theorem \ref{thm:Lyapunov}, Theorem \ref{thm:convergence_ar}. In Section \ref{sect:4}, We depict some numerical results for convex and strongly convex functions that outperform some conventional algorithms. In Section \ref{sect:5}, we conclude this paper and discuss some promising future research directions .

\section{Non-autonomous Hamiltonian Systems}\label{sect:2}

We consider the following non-autonomous \textit{Hamiltonian} system:
\begin{equation}\label{eq:H_wl}
    \begin{cases}
    \Dot{x}&=\partial_y\mathcal{H}=k(t)y\\
    \Dot{y}&=-\partial_x\mathcal{H}=-u(t)\nabla f(x)\\
\end{cases}
\end{equation}
where $\cdot=\frac{d}{dt}$, $t$ represents the time and $\mathcal{H}(x,y,t)=\frac{1}{2}k(t)||y||^2+u(t)f(x)$ is the \textit{Hamiltonian}. $K(y,t)=\frac{1}{2}k(t)||y||^2$ and $U(x,t)=u(t)f(x)$ are traditionally called the kinetic and potential energy, respectively. We assume that the time-dependent coefficients $k(t)$ and $u(t)$ are both positive when $t$ is positive. The vectors $x,\ y\in\mathbb{R}^d$ are called the position and momentum vectors, respectively. $y$ is said to be the conjugate variable of $x$. When it comes to construct numerical integrators to solve (\ref{eq:H_wl}), we usually make it autonomous by introducing a new variable $\epsilon\in\mathbb{R}$ and consider the extended \textit{Hamiltonian} \cite{feng2010symplectic}:
\begin{equation}\label{expression:H_extended}
    \Tilde{\mathcal{H}}(x,t,y,\epsilon)=\mathcal{H}(x,y,t)+\epsilon.
\end{equation}
The readers can refer to \cite{blanes2012splitting} and references therein for more proposals for making non-autonomous Hamiltonian autonomous. X=$
\begin{pmatrix}
  x\\
  t
\end{pmatrix}
$, $Y=\begin{pmatrix}
  y\\
  \epsilon
\end{pmatrix}\in\mathbb{R}^{d+1}$ are the extended position and momentum variables. $y$ is conjugate to $x$, while $\epsilon$ is conjugate to $t$. Let $\tau$ denote the time of the extended space and $'=\frac{d}{d\tau}$. Then the extended \textit{Hamiltonian} system is
\begin{equation}\label{eq:H_extended}
    \begin{cases}
     x^{'}=\partial_y\Tilde{\mathcal{H}}=\partial_y\mathcal{H}=k(t)y\\
     y^{'}=-\partial_x\Tilde{\mathcal{H}}=-\partial_x\mathcal{H}=-u(t)\nabla f(x)\\
     t^{'}=\partial_\epsilon\Tilde{\mathcal{H}}=1\\
     \epsilon^{'}=-\partial_t\Tilde{\mathcal{H}}=-\partial_t\mathcal{H}=-\frac{1}{2}\Dot{k}(t)||y||^2-\Dot{u}(t)f(x)
    \end{cases}
\end{equation}
whose projection onto the $(x,y)$-space is equal to  (\ref{eq:H_wl}) up to a time translation. And we can find that (\ref{eq:H_extended}) is an autonomous system. Therefore, the flow of (\ref{eq:H_extended}) preserves the following closed 2-form \cite{arnol2013mathematical}:
\begin{equation}
\Tilde{\omega}=dY\wedge dX=\sum_{i=1}^d dy_i\wedge dx_i+d\epsilon\wedge dt.
\end{equation}
We aim to construct numerical integrators to approximate the exact flow which also preserve $\Tilde{\omega}$. $\Tilde{\omega}$ is called a symplectic structure on $\mathbb{R}^{2(d+1)}$. ($\mathbb{R}^{2(d+1)}$,$\Tilde{\omega}$) is a symplectic manifold.
\begin{definition}
    A differentiable map $g:(\mathbb{R}^{2(d+1)},\Tilde{\omega})\to (\mathbb{R}^{2(d+1)},\Tilde{\omega})$ is said to be symplectic if it preserves the symplectic structure on $\mathbb{R}^{2(d+1)}$ i.e. $g^*\Tilde{\omega}=\Tilde{\omega}$, where $g^*$ is the pull-back of $g$.
\end{definition}
Hence the composition of symplectic maps is still symplectic. We can apply splitting methods to construct symplectic integrators \cite{hairer2006geometric}. The extended \textit{Hamiltonian} (\ref{expression:H_extended}) can be split into two parts:
\begin{equation}
    \Tilde{\mathcal{H}}(x,t,y,\epsilon)=\Tilde{\mathcal{H}}_1(x,t,y,\epsilon)+\Tilde{\mathcal{H}}_2(x,t,y,\epsilon)
\end{equation}
where $\Tilde{\mathcal{H}}_1=\frac{1}{2}k(t)||y||^2+\epsilon$ and $\Tilde{\mathcal{H}}_2=u(t)f(x)$. Then we have the following \textit{Hamiltonian} subsystems with respect to $\Tilde{\mathcal{H}}_1$:
\begin{equation}\label{eq:H_split_1}
    \begin{cases}
     x^{'}=\partial_y\Tilde{\mathcal{H}}_1=k(t)y\\
     y^{'}=-\partial_x\Tilde{\mathcal{H}}_1=0\\
     t^{'}=\partial_\epsilon\Tilde{\mathcal{H}}_1=1\\
     \epsilon^{'}=-\partial_t\Tilde{\mathcal{H}}_1=-\frac{1}{2}\Dot{k}(t)||y||^2
    \end{cases}
\end{equation}
and with respect to $\Tilde{\mathcal{H}}_2$:
\begin{equation}\label{eq:H_split_2}
    \begin{cases}
     x^{'}=\partial_y\Tilde{\mathcal{H}}_2=0\\
     y^{'}=-\partial_x\Tilde{\mathcal{H}}_2=-u(t)\nabla f(x)\\
     t^{'}=\partial_\epsilon\Tilde{\mathcal{H}}_2=0\\
     \epsilon^{'}=-\partial_t\Tilde{\mathcal{H}}_2=-\Dot{u}(t)f(x)
    \end{cases}
\end{equation}
 which can be either exactly computed or accurately approximated. Let $\Tilde{\varphi}_1^{h}$ and $\Tilde{\varphi}_2^{h}$ denote the exact $h$-flows of systems (\ref{eq:H_split_1}) and (\ref{eq:H_split_2}), respectively. Then we can integrate the above systems and calculate that
 \begin{equation}\label{expression:flow_1}
     \Tilde{\varphi}_1^{h}(x_0,y_0,t_0,\epsilon_0)=(x_0+(\int_{t_0}^{t_0+h} k(\tau)d\tau) y_0, y_0, t_0+h, \epsilon_0-\frac{1}{2}(k(t_0+h)-k(t_0))||y_0||^2)
 \end{equation}
 \begin{equation}\label{expression:flow_2}
     \Tilde{\varphi}_2^{h}(x_0,y_0,t_0,\epsilon_0)=(x_0,y_0-hu(t_0)\nabla f(x_0),t_0,\epsilon_0-h\Dot{u}(t_0)f(x_0))
 \end{equation}
and the exact flow of (\ref{eq:H_extended}) can be approximated by a composition of the form 
\begin{equation}
    \Tilde{\Psi}^h=\Tilde{\varphi}_1^{a_kh}\circ \Tilde{\varphi}_2^{b_kh}\circ\cdots\circ\Tilde{\varphi}_1^{a_1h}\circ\Tilde{\varphi}_2^{b_1h},
\end{equation}
where the coefficients $a_i,\ b_i,\ i=1,\ldots,k$ can be chosen to approximate the exact flow to some order. We consider the most commonly used and simplest one
\begin{equation}\label{alg:SE}
    \Tilde{\Psi}^h=\Tilde{\varphi}_1^{h}\circ \Tilde{\varphi}_2^{h},
\end{equation}
which is a first-order method. We can notice that there is an integral in (\ref{expression:flow_1}), which may not be exactly computed. If we apply a numerical integrator to the integral , it can still be symplectic from the following Lemma.
\begin{lemma}
    If $N^h(t_0;k)$ is a numerical integrator of the integral $\int_{t_0}^{t_0+h}k(\tau)d\tau$, and is differentiable with respect to $t_0$, then 
    \begin{equation}
        \Tilde{\psi}^h_1:(x_0,y_0,t_0,\epsilon_0)\to(x_1,y_1,t_1,\epsilon_1)=(x_0+N^h(t_0;k)y_0,y_0,t_0+h,\epsilon_0-\frac{1}{2}\Dot{N}^h(t_0;k)||y_0||^2)
    \end{equation}
    is symplectic.
\end{lemma}
\begin{proof}
From the above definition of $ \Tilde{\psi}^h_1$, we can directly calculate that
    $dx_1=dx_0+N^h(t_0)dy_0+\Dot{N}^h(t_0)y_0dt_0$, $dy_1=dy_0$, $dt_1=dt_0$ and $d\epsilon_1=d\epsilon_0-\Dot{N}^h(t_0)y_0dy_0-\frac{1}{2}\Ddot{N}^h(t_0)||y_0||^2dt_0$. Then it follows that
    \begin{align*}
        (\Tilde{\psi}^h_1)^*\Tilde{\omega}&=dy_1\wedge dx_1+d\epsilon_1\wedge dt_1\\
        &=dy_0\wedge dx_0+\Dot{N}^h(t_0)y_0dy_0\wedge dt_0+d\epsilon_0\wedge dt_0-\Dot{N}^h(t_0)y_0dy_0\wedge dt_0\\
        &=dy_0\wedge dx_0+d\epsilon_0\wedge dt_0=\Tilde{\omega}.
    \end{align*}
    Therefore,  $ \Tilde{\psi}^h_1$ is symplectic.
    \end{proof}
For simplicity, we just apply the rectangle integrator $N^h(t_0;k)=hk(t_0)$, then the method (\ref{alg:SE}) is modified as follows:
\begin{equation}
   \Tilde{\Psi}^h=\Tilde{\psi}_1^{h}\circ \Tilde{\varphi}_2^{h}
\end{equation}
i.e.
\begin{equation}
    \Tilde{\Psi}^h: (x_0,y_0,t_0,\epsilon_0)\to(x_1,y_1,t_1,\epsilon_1)=(x_0+hk(t_0)y_1,y_0-hu(t_0)\nabla f(x_0),t_0+h,\epsilon_0-\frac{1}{2}h\Dot{k}(t_0)||y_1||^2-h\Dot{u}(t_0)f(x_0))
\end{equation}
which is still a symplectic method. Projected onto the $(x,y,t)$-space, we have the $(n+1)_{th}$ update as follows:
\begin{equation}\label{alg:SE_xyuk}
\begin{cases}
    y_{n+1}=y_{n}-h_nu(t_n)\nabla f(x_n)\\
    x_{n+1}=x_n+h_nk(t_n)y_{n+1}\\
    t_{n+1}=t_n+h_n
\end{cases}
\end{equation}
and $\{(x_n,y_n)\}_{n=0}^\infty$ is called a numerical trajectory of the scheme with step-size $\{h_n\}_{n=0}^\infty$. We can linearize it when we consider the Taylor expansion of the vector field $\begin{pmatrix}
    \nabla f(x)\\y
\end{pmatrix}$ near its critical point $\begin{pmatrix}
    x^*\\0
\end{pmatrix}$. Without loss of generality, we assume $x^*=0$, otherwise, we apply the translation $x-x^*$. Therefore, $\nabla f(x)=\nabla^2f(0)x+O(||x||^2)$. Substituting it into (\ref{alg:SE_xyuk}), we can get its linearized form if we omit the second-order terms.
\begin{definition}
    A scheme is called (linearly) stable with step-size $\{h_n\}_{n=0}^\infty$ if the corresponding trajectory $\{(x_n,y_n)\}_{n=0}^\infty$ of its linearized form remains bounded .
\end{definition}
Now we can prove Theorem \ref{thm1}:
\begin{proof}[Proof of Theorem \ref{thm1}]
Assume $x^*=0$, then $f(x)=f(x^*)+\frac{1}{2}x^T\nabla^2 f(x^*) x+O(||x||^3)$, and $\nabla^2 f(x^*)$ can be diagonalized using an orthogonal transformation $P$, which means that $P\nabla^2 f(x^*)P^T=\Lambda$ and $\Lambda=diag\{\lambda_1,...,\lambda_d\}$. Since $f(x)$ is convex and $\nabla f(x)$ is Lipschitz continuous with Lipschitz constant $L$, we have $0\leq\lambda_i\leq L$.Then we can write (\ref{alg:SE_xyuk}) as $\begin{pmatrix}
    x_{n+1}\\ y_{n+1}
\end{pmatrix}=S(t_n;h_n)\begin{pmatrix}
    x_n\\y_n
\end{pmatrix}$.
The stability matrix $S(t_n;h_n)$ is 
\begin{equation}
    S(t_n;h_n)=\begin{pmatrix}
        I_d-h_n^2k(t_n)u(t_n)\nabla^2 f(x^*) & h_nk(t_n)I_d\\
        -h_nu(t_n)\nabla^2 f(x^*) & I_d
    \end{pmatrix}=\begin{pmatrix}
        P & 0\\ 0 & P
    \end{pmatrix}\begin{pmatrix}
        I_d-h_n^2k(t_n)u(t_n)\Lambda & h_nk(t_n)I_d\\
        -h_nu(t_n)\Lambda & I_d
    \end{pmatrix}\begin{pmatrix}
        P^T & 0\\ 0 & P^T
    \end{pmatrix}
\end{equation}
where $I_d$ is the $d\times d$ identity matrix. Therefore, $S(t_n;h_n)$ can be congruent to a block-diagonal matrix consisting of $d$ $2\times 2$ matrices:
\begin{equation}
    S_i(t_n;h_n)=\begin{pmatrix}
        1-h_n^2k(t_n)u(t_n)\lambda_i & h_nk(t_n)\\
        -h_nu(t_n)\lambda_i & 1
    \end{pmatrix}, i=1,\ldots,d.
\end{equation}
Let $\sigma(A)$ denote the set of all eigenvalues of matrix $A$ and $\rho(A)=max\{|\lambda|:\lambda\in\sigma(A)\}$ is the spectral radius of $A$. The scheme (\ref{alg:SE_xyuk}) is stable if $\rho(S(t_n;h_n))\leq 1$ and $S(t_n;h_n)$ can be diagonalized. Then we have $\sigma(S(t_n;h_n))=\bigcup_{i=1}^d\sigma(S_i(t_n;h_n))$ and $S(t_n;h_n)$ can be diagonalized if all $S_i(t_n;h_n)$ can be diagonalized. When $S_i(t_n;h_n)$ has two distinct eigenvalues, it can be diagonalized. Since $\det S_i(t_n;h_n)=1$, if it has two distinct real eigenvalues, then the eigenvalues must be reciprocals of each other, hence $\rho(S_i(t_n;h_n))>1$. It follows that $\rho(S_i(t_n;h_n))\leq1$ and $\#\sigma(S_i(t_n;h_n))=2$ is equivalent to $|tr(S_i(t_n;h_n))|<2$ i.e. $h_n^2k(t_n)u(t_n)<\frac{4}{\lambda_i}\leq\frac{4}{L}$. We can get the scheme (\ref{alg:SE_ab}) from (\ref{alg:SE_xyuk}) when $k(t)=\frac{1}{a(t)}$ and $u(t)=b(t)$. Therefore, (\ref{conditions:stepsize}) can be derived from the above process.
\end{proof}
\begin{proof}[Proof of Corollary \ref{cor1}]
Considering the polynomial subfamily of the \textit{Bregman} \textit{Lagrangian} in the \textit{Euclidean setting}, we have $a(t)=e^{\gamma(t)-\alpha(t)}=\frac{1}{p}t^{p+1}$ and $b(t)=e^{\alpha(t)+\beta(t)+\gamma(t)}=Cpt^{2p-1}$ from (\ref{conditions:polynomial}). Therefore, $\frac{b(t)}{a(t)}=Cp^2t^{p-2}$, and we can get (\ref{conditions:stepsize_poly}) by substituting it into the (\ref{conditions:stepsize}) in Theorem \ref{thm1}.
\end{proof}
We illustrate with a simple example. Let $f(x)=\frac{1}{2}x^TM^{-1}x$, where $M=(m_{ij})\in\mathbb{R}^{d\times d}$, $m_{ij}=0.9^{|i-j|}$, $x\in \mathbb{R}^d$ and $d=50$. We apply the scheme (\ref{alg:SE_xyuk}) to the polynomial subfamily of the \textit{Bregman Lagrangians} i.e. $u(t)=Cpt^{2p-1}$, $k(t)=pt^{-(p+1)}$ with a fixed step-size and varied step-sizes satisfying the stable condition (\ref{conditions:stepsize_poly}). In particular, we choose $C=\frac{1}{Lp^2}$, $L$ is the Lipschitz constant of $\nabla f$, then $h_n=\sqrt{\frac{1}{t_n^{p-2}}}$ satisfies the stable condition from Corollary \ref{cor1}.
\begin{figure}[htb]
    \centering
\subfigure[
]{ \label{fig:Breg_h}
\begin{minipage}[b]{0.3\textwidth}
    \includegraphics[width=\textwidth]{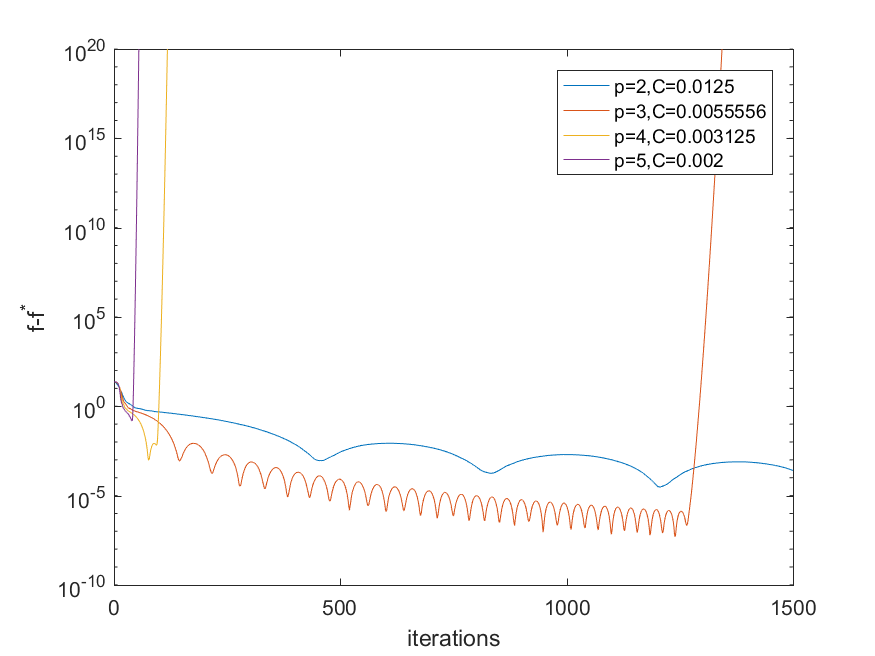}
  \end{minipage}
   }
  \subfigure[
  ]{\label{fig:Breg_itrations}
    \begin{minipage}[b]{0.3\textwidth}
    \includegraphics[width=\textwidth]{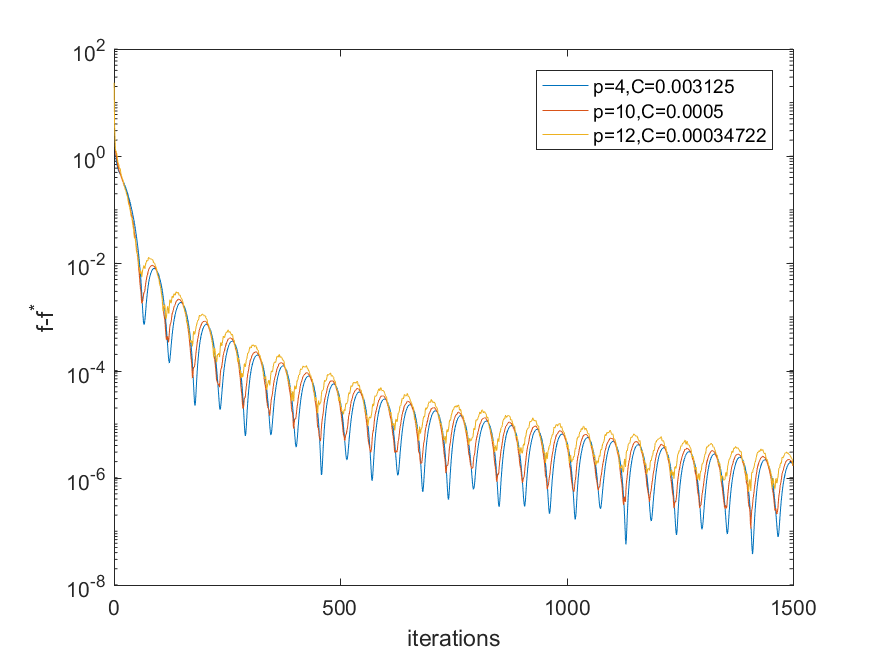}
  \end{minipage}
       }
        \subfigure[
  ]{\label{fig:Breg_time}
    \begin{minipage}[b]{0.3\textwidth}
    \includegraphics[width=\textwidth]{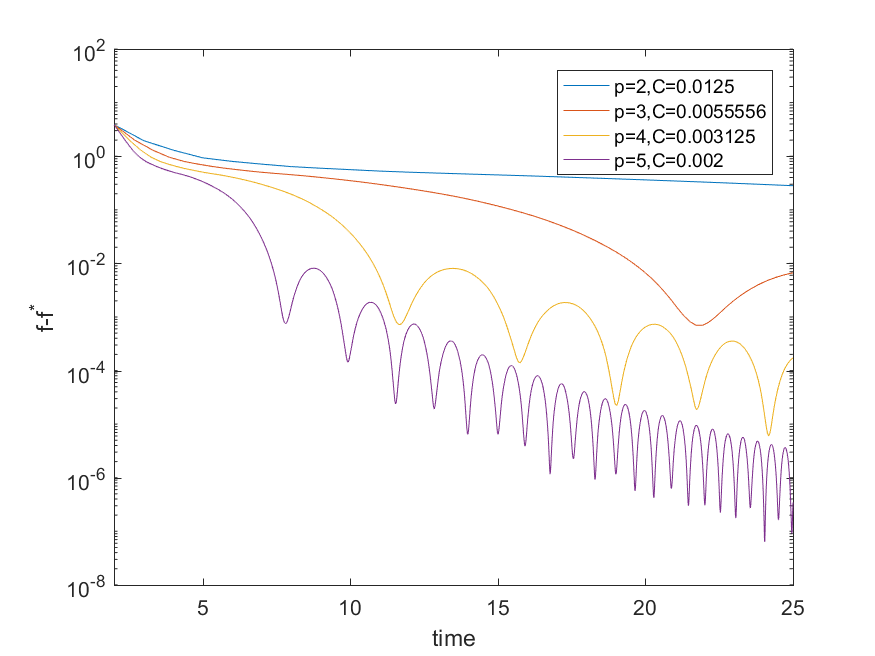}
  \end{minipage}
       }
\caption{Scheme (\ref{alg:SE_xyuk}) with fixed and varied step-size}
\end{figure}\par
 As shown in \ref{fig:Breg_h}, (\ref{alg:SE_xyuk}) with a fixed step-size will become unbounded when the iterations become large enough for $p
 \geq 3$, while it remains bounded for $p=2$. Besides, with the same step-size, the numerical trajectory requires fewer iterations to diverge as $p$ increases. In contrast, if we take step-size in the stable region which is updated with each iteration, the numerical solutions are all bounded for $p\geq 2$ as depicted in \ref{fig:Breg_itrations} and \ref{fig:Breg_time}. We take the iteration $n$ and the time $t_n$ as x-axis in \ref{fig:Breg_itrations} and \ref{fig:Breg_time}, respectively. From the iteration point of view, (\ref{alg:SE_xyuk}) exhibits similar oscillatory period, but smaller amplitude as $p$ gets larger. However, if we view the numerical solutions from the time perspective, they are consistent with the continuous case. As $p$ increases, the convergence rate becomes faster and the oscillations become more frequent.

\section{Lyapunov Analysis}\label{sect:3}
We consider the following $(\alpha,r)$-damped system:
\begin{equation}\label{eqn:E-L_ar2}
    \Ddot{x}+\frac{r}{t^\alpha}\Dot{x}+\nabla f(x)=0
\end{equation}
where $r>0$ is and $0\leq \alpha \leq 1$.
First of all, we need to verify the global existence and uniqueness of solution to the equation (\ref{eqn:E-L_ar2}).
\begin{lemma}\label{lem:dissipative}
    The system (\ref{eqn:E-L_ar2}) is a dissipative system i.e. its energy
   \begin{equation}\label{energy}
       \mathcal{E}(t)=\frac{1}{2}||\Dot{x}(t)||^2+(f(x(t))-f^*)
   \end{equation}
   is decreasing along any trajectory of the system. 
\end{lemma}
\begin{proof}
    Differentiating (\ref{energy}) with respect to $t$ and substituting equation (\ref{eqn:E-L_ar2}), then
    \begin{align*}
    \frac{d}{dt}\mathcal{E}(t)&=<\Ddot{x},\Dot{x}>+<\nabla f(x),\Dot{x}>\\
    &=-\frac{r}{t^\alpha}||\Dot{x}||^2\leq 0.
    \end{align*}
\end{proof}
\begin{theorem}\label{thm:existence}
    If $f$ is twice continuously differentiable and $\nabla f$ is globally Lipschitz continuous, then for any $x_0\in \mathbb{R}^d$ and $t_0>0$, the ODE (\ref{eqn:E-L_ar2}) has a unique global solution $x(t):[t_0,\infty)\to\mathbb{R}^d$ satisfying the initial conditions $x(t_0)=x_0$ and $\Dot{x}(t_0)=0$. 
\end{theorem}
\begin{proof}
    We rewrite the second-order ODE (\ref{eqn:E-L_ar2}) as the following first-order ODEs:
    \begin{equation}
        \Dot{x}=v,\ \Dot{v}=-\frac{r}{t^\alpha}v-\nabla f(x).
    \end{equation}
   Then the vector field $\begin{pmatrix}
        v\\-\frac{r}{t^\alpha}v-\nabla f(x))
    \end{pmatrix}$ is globally Lipschitz continuous with respect to $(x,v)$ uniformly in $t$,  since $\frac{r}{t^\alpha}$ is bounded when $t\geq t_0>0$ and $\nabla f$ is globally Lipschitz continuous. 
    Therefore, the ODE (\ref{eqn:E-L_ar2}) has a unique global solution satisfying the initial conditions $x(t_0)=x_0$ and $\Dot{x}(t_0)=0$ referring to Theorem 3 of Chapter 3.1 in \cite{perko2013differential}. 
\end{proof}
\begin{definition}
    $f:\mathbb{R}^d\to \mathbb{R}$ is called a \textit{convex} function if $\forall x,\ y \in \mathbb{R}^d$
    \begin{equation}
        f(y)-f(x)-<\nabla f(x),y-x>\geq 0.
    \end{equation}
\end{definition}
\begin{definition}
    $f:\mathbb{R}^d\to \mathbb{R}$ is called a $\mu$-\textit{strongly convex} function if $\forall x,\ y \in \mathbb{R}^d$
    \begin{equation}
        f(y)-f(x)-<\nabla f(x),y-x>\geq 
        \frac{\mu}{2}||y-x||^2.
    \end{equation}
\end{definition}
Motivated by \cite{li2024linear}, we assume the Lyapunov function as follows:
\begin{equation}\label{Lyapunov2}
    L(t)=\frac{1}{2}||A(t)\Dot{x}(t)+B(t)(x(t)-x^*)||^2+\frac{1}{2}C(t)||\Dot{x}(t)||^2+D(t)(f(x(t))-f^*).
\end{equation}
We aim to find $A,\ B,\ C,\ D:[t_0,\infty)\to [0,\infty)$ such that $\frac{d}{dt}L(t)\leq 0$ along the solution of the equation (\ref{eqn:E-L_ar}) and $D(t)$ is as large as enough.\\
Differentiating (\ref{Lyapunov2}) with respect to $t$ and substituting equation (\ref{eqn:E-L_ar2}), then 
\begin{align}
    \frac{d}{dt}L&=[\frac{1}{2}(2A\Dot{A}+\Dot{C})-\Dot{\xi}(A^2+C)+AB]||\Dot{x}||^2\tag{\Rmnum{1}}\label{1}\\
    &+B\Dot{B}||x-x^*||^2\tag{\Rmnum{2}$_{1}$}\label{2-1}\\
    &+(\Dot{D}-AB)(f(x)-f^*)\tag{\Rmnum{2}$_{2}$}\label{2-2}\\
    &-AB[f^*-f(x)-<\nabla f(x),x^*-x>]\tag{\Rmnum{2}$_{3}$}\label{2-3}\\
    &+[D-A^2-C]<\nabla f(x),\Dot{x}>\tag{\Rmnum{3}}\label{3}\\
    &+[A\Dot{B}+\Dot{A}B-\Dot{\xi}AB+B^2]<\Dot{x},x-x^*>\tag{\Rmnum{4}}\label{4}
\end{align}
where $\Dot{\xi}=\frac{r}{t^\alpha}$.
\begin{lemma} \label{condition1}
If the following conditions are satisfied
    \begin{gather}
        \frac{1}{2}\Dot{D}-\Dot{\xi}D+AB\leq 0,\tag{\Rmnum{1}$^{'}$}\label{1'}\\
        B\Dot{B}\leq 0, \tag{\Rmnum{2}$_{1}^{'}$}\label{2-1'}\\
        \Dot{D}-AB\leq 0, \tag{\Rmnum{2}$_{2}^{'}$}\label{2-2'}\\
        D-A^2-C=0, \tag{\Rmnum{3}$^{'}$}\label{3'}\\
        \frac{d}{dt}(AB)-\Dot{\xi}AB+B^2=0\tag{\Rmnum{4}$^{'}$}\label{4'},
    \end{gather}
    then $\frac{d}{dt}L(t)\leq 0$.
\end{lemma}

\begin{proof}
    $f$ is convex, then $f(y)-f(x)-<\nabla f(x),y-x>\geq 0$ from the definition. Let $y=x^*$, we have $f^*-f(x)-<\nabla f(x),x^*-x>\geq 0$. Hence (\ref{2-3}) is always negative when $AB$ is positive. (\ref{3'}) and (\ref{4'}) are derived from (\ref{3}) and (\ref{4}), since $<\nabla f(x),\Dot{x}>$ and $<\Dot{x},x-x^*>$ can be positive and negative, their coefficients must be zero. (\ref{1}), (\ref{2-1}) and (\ref{2-2}) are negative if their coefficients are negative. Thus we can induced (\ref{1'}), (\ref{2-1'}) and (\ref{2-2'}) from this and (\ref{3'}).
\end{proof}
\begin{assumption}\label{a1}
    $D(t)=e^{\beta\xi(t)}$, $\beta\geq 0$.
\end{assumption}
\begin{assumption}\label{a2}
    $\Dot{D}=AB$.
\end{assumption}
\begin{lemma}\label{conditions:beta}
    Under the above assumptions, if 
    \begin{equation}
        \beta\leq \begin{cases}
            \frac{1}{2}\quad,\ &0\leq \alpha<1\\
            \min\{\frac{2}{3},\frac{1+r}{2r}\},\ &\alpha=1,
        \end{cases}
    \end{equation}
     then (\ref{1'}) is true and we can solve $A,\ B,\ C$ satisfying (\ref{2-2'}), (\ref{3'}) and (\ref{4'}) and with the expressions (\ref{expression:At}), (\ref{expression:Bt}), (\ref{expression:Ct}) stated in the Theorem \ref{thm:Lyapunov} .
\end{lemma}
\begin{proof}
    Substituting Assumption \ref{a2} into (\ref{1'}), then
    \begin{equation}
        \frac{3}{2}\Dot{D}-\Dot{\xi}D\leq 0.
    \end{equation}
   Applying the Gronwall's inequality, we get
   \begin{equation}
       D(t)\leq e^{\frac{2}{3}\xi(t)},
   \end{equation}
   then combining with Assumption \ref{a1}, we have $\beta\leq\frac{2}{3}$. Substituting Assumption \ref{a2} into (\ref{4'}), then 
   \begin{equation}\label{B}
       B^2=\Dot{\xi}\Dot{D}-\Ddot{D}
   \end{equation}
   By Assumption \ref{a1}, $D(t)=e^{\beta\xi(t)}$, therefore
   \begin{equation}\label{Ddot}
       \Dot{D}=\beta\Dot{\xi}e^{\beta\xi(t)},\ \Ddot{D}=[\beta^2\Dot{\xi}^2+\beta\Ddot{\xi}]e^{\beta\xi(t)},
   \end{equation}
   where $\Dot{\xi}=rt^{-\alpha},\ \Ddot{\xi}=-r\alpha t^{-\alpha-1}$. Substituting (\ref{Ddot}) into (\ref{B}), we can get the expression of $B^2$:
   \begin{align}
       B^2&=[\beta(1-\beta)\Dot{\xi}^2-\beta\Ddot{\xi}]e^{\beta\xi(t)}\label{B2}\\
       &=\beta r[(1-\beta)rt^{1-\alpha}+\alpha]t^{-\alpha-1}e^{\beta\xi(t)}\geq 0,
   \end{align}
   then we arrive at the expression of $B$ (\ref{expression:Bt}). And by Assumption \ref{a2}, we can infer that 
       \begin{equation}
           A=\frac{\Dot{D}}{B}=\sqrt{\frac{\beta rt^{1-\alpha}}{(1-\beta)rt^{1-\alpha}+\alpha}e^{\beta\xi(t)}}\geq 0,
       \end{equation}
       which is exactly (\ref{expression:At}). Now we obtain $A,\ B,\ D$ satisfying (\ref{1'}), (\ref{2-2'}) and (\ref{4'}). It remains to solve $C\geq 0$ from (\ref{3'}):
       \begin{equation}
           C=D-A^2=[1-\frac{\beta rt^{1-\alpha}}{(1-\beta)rt^{1-\alpha}+\alpha}]e^{\beta\xi(t)},
       \end{equation}
       which is (\ref{expression:Ct}). And $C\geq 0$ if $(1-\beta)rt^{1-\alpha}+\alpha\geq\beta r t^{1-\alpha}$. There are two cases as follows:\\
       Case 1: $0\leq\alpha<1$, then it is equivalent to $(1-\beta)r\geq \beta r $ i.e. $\beta\leq\frac{1}{2}$;\\
       Case 2: $\alpha=1$, then we have $(1-\beta)r+1\geq \beta r$ i.e. $\beta\leq \frac{1+r}{2r}$.\\
       Therefore, we can get $A,\ B,\ C,\ D$ satisfying (\ref{1'}), (\ref{2-2'}), (\ref{3'}) and (\ref{4'}) if $\beta\leq \min\{\frac{2}{3},\frac{1}{2}\}=\frac{1}{2}$ when $0\leq\alpha<1$, or $\beta\leq \min\{\frac{2}{3},\frac{1+r}{2r}\}$ when $\alpha=1$.
\end{proof}
Now it remains to verify (\ref{2-1'}), however, we have the following Lemma:
\begin{lemma} \label{lemma:why strongly convex}
Under the above assumptions, we have\\
    (i) If $0\leq \alpha<1$, (\ref{2-1'}) and (\ref{4'}) cannot be simultaneously true.\\
    (ii) If $\alpha=1$, (\ref{2-1'}) and (\ref{4'}) can be both true when $\beta r\leq 2$.
\end{lemma}
\begin{proof}
    We assume that (\ref{4'}) is true. From (\ref{4'}) and Assumption \ref{a2}, as in the proof of Lemma \ref{conditions:beta}, we get (\ref{B2}). Differentiating both sides of it, then
    \begin{align}
        2B\Dot{B}&=[\beta^2(1-\beta)\Dot{\xi}^3+\beta(2-3\beta)\Dot{\xi}\Ddot{\xi}-\beta\overset{\dots}{\xi}]e^{\beta\xi}\\
        &=\beta r[\beta(1-\beta)r^2t^{2(1-\alpha)}-(2-3\beta)r\alpha t^{1-\alpha}-\alpha(1+\alpha)]t^{-\alpha-2}e^{\beta\xi}\label{Bdot}.
    \end{align}
    We can infer that $B\Dot{B}>0$ when $0\leq\alpha<1$, $t>t_{\alpha,r,\beta}$, and $t_{\alpha,r,\beta}^{1-\alpha}$ is the positive solution of
    \begin{equation}
        \beta(1-\beta)r^2\lambda^2-(2-3\beta)r\alpha \lambda-\alpha(1+\alpha)=0.
    \end{equation}
    Therefore, (\ref{2-1'}) cannot be satisfied if $0\leq\alpha<1$ and $t$ is large. When $\alpha=1$, $B^2=\beta r[(1-\beta)r+1]t^{\beta r-2}$ is decreasing, thus  (\ref{2-1'}) holds when $\beta r\leq 2$. 
\end{proof}
From the above Lemma, we cannot construct a Lyapunov function of the form (\ref{Lyapunov2}) satisfying Lemma \ref{condition1} for convex function $f$ when $\alpha<1$ , which motivates us to consider $\mu$-strongly convex function.
\begin{lemma}\label{conditions:mu}
    If $f$ is a $\mu$-strongly convex function, then
   \begin{center}
        (\Rmnum{2})=(\ref{2-1})+(\ref{2-2})+(\ref{2-3})$\leq 0$
   \end{center}
   when $t\geq T$ and $T=T(\mu,\alpha,r,
    \beta)$ is either the largest root of the polynomial 
    \begin{equation}\label{expression:polynomial_P}
        P(t;\mu,\alpha,r,\beta)=\mu t^2-(1-\beta)\beta r^2t^{2(1-\alpha)}+(2-3\beta)r\alpha t^{1-\alpha}+\alpha(\alpha+1),
    \end{equation}
     or zero if $P(t;\mu,\alpha,r,\beta)$ has no root. In particular, when $\alpha=0$, $T=0$ and $r$ should satisfy that $r\leq\sqrt{\frac{\mu}{\beta(1-\beta)}}$.
\end{lemma}
\begin{proof}
    For a $\mu$-strongly convex function $f$, we have
    \begin{equation}
        f^*-f(x)-<\nabla f(x),x^*-x>\geq \frac{\mu}{2}||x-x^*||^2
    \end{equation}
    Under the Assumption \ref{a2}, we have
    \begin{align}
        (\Rmnum{2})&=B\Dot{B}||x-x^*||^2-\Dot{D}[f^*-f(x)-<\nabla f(x),x^*-x>]\\
        &\leq \frac{1}{2}[2B\Dot{B}-\mu\Dot{D}]||x-x^*||^2\\
        &=-\frac{1}{2}\beta r [\mu t^2-(1-\beta)\beta r^2t^{2(1-\alpha)}+(2-3\beta)r\alpha t^{1-\alpha}+\alpha(\alpha+1)]t^{-\alpha-2}e^{\beta\xi}\\
        &=-\frac{1}{2}\beta r P(t;\mu,\alpha,r,\beta)t^{-\alpha-2}e^{\beta\xi}
    \end{align}
    where the second equation is derived from (\ref{Bdot}) and $\Dot{D}=\beta r t^{-\alpha}e^{\beta\xi}$. Therefore, $(\Rmnum{2})\leq 0$ if the polynomial $P(t;\mu,\alpha,r,\beta)\geq 0$. When $\alpha>0$, $P(t;\mu,\alpha,r,\beta)> 0$ as $t\to \infty$, then we only need that $t$ is larger than the largest root of $P$. For $\alpha=0$, $P=(\mu-(1-\beta)\beta r^2)t^2$, then $P\geq 0$ only when $r\leq \sqrt{\frac{\mu}{\beta(1-\beta)}}$.
\end{proof}
We can get the exact expression of $T$ when $\alpha=1$:
\begin{lemma}\label{expresion:T_a_1}
    When $\alpha=1$, we have
    \begin{equation}
        T(\mu,r,1,\beta)=
        \begin{cases}
            0,\ \beta r\leq 2\\
            \sqrt{\frac{(\beta r-2)(r+1-\beta r)}{\mu}},\ \beta r >2.
        \end{cases}
    \end{equation}
\end{lemma}
\begin{proof}
    When $\alpha=1$, $P(t;\mu,r,1,\beta)=\mu t^2-(1-\beta)\beta r^2+(2-3\beta)r+2$. \\
    Let $k=\beta r$, $\delta=-(1-\beta)\beta r^2+(2-3\beta)r+2=(k-2)(k-r-1)$, the roots of $\delta(k)$ are $2$ and $r+1$. From Lemma \ref{conditions:beta}, we have $k\leq \min\{\frac{2}{3}r,\frac{r+1}{2}\}$. If $r\leq 3$, then $k\leq \frac{2}{3}r\leq min\{2,r+1\}$. Therefore, $\delta\geq 0$ when $r\leq 3$. If $r>3$, then $k\leq \frac{r+1}{2}<r+1$. We can deduce that $\delta\leq 0$ when $k\geq 2$ and $\delta> 0$ when $k<2$. \\
    In conclusion, if $k< 2$, then $\delta> 0$, $P$ has no root. If $k=2$, then $\delta=0$, the only root of $P$ is zero. Hence $T=0$ when $\beta r=k\leq 2$. If $\beta r=k>2$, then $\delta<0$, the roots of $P$ are $\pm\sqrt{-\delta}$, $T=\sqrt{-\delta}=\sqrt{\frac{(\beta r-2)(r+1-\beta r)}{\mu}}$.
\end{proof}
\begin{remark}\label{rmk:expression:T}
    For $0<\alpha<1$, we can get an approximate expression of $T$. Consider the main terms of $P(t;
    \mu,\alpha,r,\beta)$ i.e. $\mu t^2-(1-\beta)\beta r^2t^{2(1-\alpha)}$ with $\mu<<1$, then $$\mu T^2\sim (1-\beta)\beta r^2T^{2(1-\alpha)}.$$
    Therefore, we can estimate that $T\sim (\frac{r^2}{\mu})^\frac{1}{2\alpha}$ which increases as $\alpha$ decreases and $r$ increases.
\end{remark}
Now we can prove our main results using the above Lemmas:
\begin{proof}[Proof of Theorem \ref{thm:Lyapunov}]
Combining Lemma \ref{conditions:beta} with Lemma \ref{conditions:mu}, we can simplify the derivative of the function $L$ in (\ref{thm:Lyapunov}) with $A,\ B,\ C,\ D $ defined in (\ref{expression:At}), (\ref{expression:Bt}), (\ref{expression:Ct}) and (\ref{expression:Dt}) as follows:
\begin{equation}
    \frac{d}{dt}L(t)=-[(1-\frac{3}{2}\beta)rt^{-\alpha}||\Dot{x}||^2+\frac{1}{2}\beta r P(t;\mu,\alpha,r,\beta)t^{-\alpha-2}||x-x^*||^2]e^{\beta\xi}\leq 0
\end{equation}
when $t\geq T$. The constraints in (\ref{conditions:beta_r}) are derived directly from $\Dot{\xi}=rt^{-\alpha}$, Lemma \ref{conditions:beta} and Lemma \ref{conditions:mu}.
\end{proof}
\begin{proof}[Proof of Theorem \ref{thm:convergence_ar}]
    (i) $\alpha=0$, we apply Theorem \ref{thm:Lyapunov} and take $\beta=\frac{1}{2},\ r=\sqrt{\frac{\mu}{\beta(1-\beta)}}=2\sqrt{\mu}$. Then $D(t)=e^{\sqrt{\mu}t}$, $B(t)=\sqrt{\mu e^{\sqrt{\mu}t}}$, $A(t)=\sqrt{e^{\sqrt{\mu}t}}$ and $C(t)=0$. Substituting into (\ref{Lyapunov2}), then $$L(t)=[\frac{1}{2}||\Dot{x}+\sqrt{\mu}(x-x^*)||^2+(f(x)-f^*)]e^{\sqrt{\mu}t}$$
    is a Lyapunov function when $t\geq T=0$ from Theorem \ref{thm:Lyapunov} and Lemma \ref{conditions:mu}. Then $e^{\sqrt{\mu}t}[f(x(t))-f^*]\leq L(t)\leq L(0)=\frac{\mu}{2}||x_0-x^*||^2+(f(x_0)-f^*)$, thus the convergence rate is
    \begin{equation}
        f(x(t))-f^*\leq[\frac{\mu}{2}||x_0-x^*||^2+f(x_0)-f^*]e^{-\sqrt{\mu}t}.
    \end{equation}
    (ii) $0<\alpha<1$, let $\beta=\frac{1}{2}$ and apply Theorem \ref{thm:Lyapunov}, we have $D(t)[f(x(t))-f^*]\leq L(t)\leq L(T)$ and $D(t)=e^{\beta \xi}=e^{\frac{r}{2(1-\alpha)}t^{1-\alpha}}$, thus the convergence rate is
    \begin{equation}
        f(x(t))-f^*\leq L(T)e^{-\frac{r}{2(1-\alpha)}t^{1-\alpha}}
    \end{equation}
    when $t\geq T$.\\
    (iii) $\alpha=1$, similar as the above cases, if $r>3$ and $\beta=\frac{r+1}{2r}$, then $D(t)=t^{\frac{r+1}{2}}$ and the convergence rate is
    \begin{equation}
        f(x(t))-f^*\leq L(T)t^{-\frac{r+1}{2}}
    \end{equation}
    since $L(t)$ is a Lyapunov function when $t\geq T$. Moreover, $T=\sqrt{\frac{(\beta r-2)(r+1-\beta r)}{\mu}}=\sqrt{\frac{(r-3)(r+1)}{4\mu}}$, from Lemma \ref{expresion:T_a_1} when $\beta r=\frac{r+1}{2} >2$.\\
    If $r\leq 3$, let $\beta=\frac{2}{3}$, then $D(t)=t^{\frac{2}{3}r}$, $B(t)=\sqrt{\frac{2r}{3}(\frac{r}{3}+1)t^{\frac{2}{3}r-2}}$, $A(t)=\sqrt{\frac{2r}{r+3}t^{\frac{2}{3}r}}$ and $C(t)=\frac{3-r}{r+3}t^{\frac{2}{3}r}$. Applying Lemma \ref{expresion:T_a_1} for $\beta r\leq 3$, then $L(t)\leq L(t_0)=\frac{1}{2}B^2(t_0)||x_0-x^*||^2+D(t_0)(f(x_0)-f^*)$ when $t\geq t_0>0$. Therefore, we can derive the following inequality:
    \begin{equation}
        f(x(t))-f^*\leq [\frac{r}{3}(\frac{r}{3}+1)t_0^{\frac{2}{3}r-2}||x_0-x^*||^2+t_0^{\frac{2}{3}r}(f(x_0)-f^*)]t^{-\frac{2}{3}r}
    \end{equation}
\end{proof}
\begin{remark}
    The cases where $r=2\sqrt{\mu},\ \alpha=0$ and $r=3,\ \alpha=1$ are indeed NAG for $\mu$-strongly convex functions and convex functions, respectively. Therefore, we build a unified framework of them. For convex functions, we can obtain at most $O(t^{-2})$ convergence rate from Lemma \ref{lemma:why strongly convex} (ii). For $\mu-$strongly convex functions and $r>3,\ \alpha=1$, \cite{NAGODE} gave a faster convergence rate $O(t^{-\frac{2}{3}r})$ via a proof by induction on $r$ and solving $B$ satisfying (\ref{4})$>0$, while we assume (\ref{4})$\equiv 0$.
\end{remark}

\section{Numerical Experiments}\label{sect:4}
For the \textit{Lagrangian }
\begin{equation}
     \mathcal{L}(x,v,t)=e^{\xi(t)}[\frac{1}{2}||v||^2-f(x)],
\end{equation}
its corresponding \textit{Hamiltonian} is 
\begin{equation}
    \mathcal{H}(x,y,t)=e^{-\xi(t)}\frac{1}{2}||y||^2+e^{\xi(t)}f(x).
\end{equation}
Therefore, the Hamiltonian equations $\Dot{x}=\partial_y \mathcal{H},\ \Dot{y}=-\partial_x \mathcal{H}$ are as follows:
\begin{equation}
    \begin{cases}
        \Dot{x}&=e^{-\xi(t)}y\\
        \Dot{y}&=-e^{\xi(t)}\nabla f(x),
    \end{cases}
\end{equation}
which are equivalent to the \textit{Euler-Lagrange} equation (\ref{eqn:E-L_ar}). We apply the symplectic Euler method introduced in Section \ref{sect:2} to integrate the system starting from $x_0 \in \mathbb{R}^d,\ y_0=0$ and $t_0>0$ with step-size $h$:
\begin{equation}\label{alg:SE_ar}
    \begin{split}
        y_{n+1}&=y_{n}-he^{\xi(t_n)}\nabla f(x_n)\\
        x_{n+1}&=x_{n}+he^{-\xi(t_n)}y_{n+1}\\
        t_{n+1}&=t_n+h
    \end{split}
\end{equation}
where $\xi(t)=\frac{r}{1-\alpha}t^{1-\alpha}$. We can take $h<\frac{2}{L}$ from Theorem \ref{thm1}. 
\begin{remark}
    We can observe that $||y_n||$ is growing very fast, thus causing some troubles such as overflow in practice. Therefore, we indeed use the equivalent scheme:
    \begin{equation}\label{alg:SE_ar_velocity}
    \begin{split}
        v_{n+1}&=e^{\xi(t_{n-1})-\xi(t_n)}v_{n}-h\nabla f(x_n)\\
        x_{n+1}&=x_{n}+hv_{n+1}\\
        t_{n+1}&=t_n+h
    \end{split}
\end{equation}
where $v_n=e^{-\xi(t_{n-1})}y_n$ represents the velocity of the $n_{th}$ iteration.
\end{remark}
We consider the following two cases: $\mu$-strongly convex functions with $\frac{\mu}{L}<<1$ and convex but not strongly convex functions as in \cite{NAGODE}.\\
\textbf{Quadratic}. $f(x)=\frac{1}{2}x^TAx+b^Tx$, where $A$ is a random positive definite matrix and $b$ is a random vector. We assume the dimension of the space $d=500$ and the eigenvalues of $A$ are between 0.001 and 1. The vector $b$ is generated as i.i.d. Gaussian random variables with mean 0 and variance 25. Therefore, $f$ is a strongly convex function with $\frac{\mu}{L}=0.001$.\\
\textbf{Log-sum-exp}. 
$f(x)=\rho \log [\sum_{i=1}^m \exp((a_i^Tx-b_i)/\rho)]$,
where $m=200,\ \rho=20$ and the dimension $d=50$. The matrix $A=(a_{ij})$ is a random matrix with i.i.d. standard Gaussian entries and $b=(b_i)$ has i.i.d. Gaussian entries with mean 0 and variance 2. $f(x)$ is a convex but not strongly convex function.\\

\begin{figure}[hbt]
\centering
\subfigure[min \text{$\frac{1}{2}x^TAx+b^Tx$}
]{ \label{fig:qudratic_a}
\begin{minipage}[b]{0.45\textwidth}
    \includegraphics[width=\textwidth]{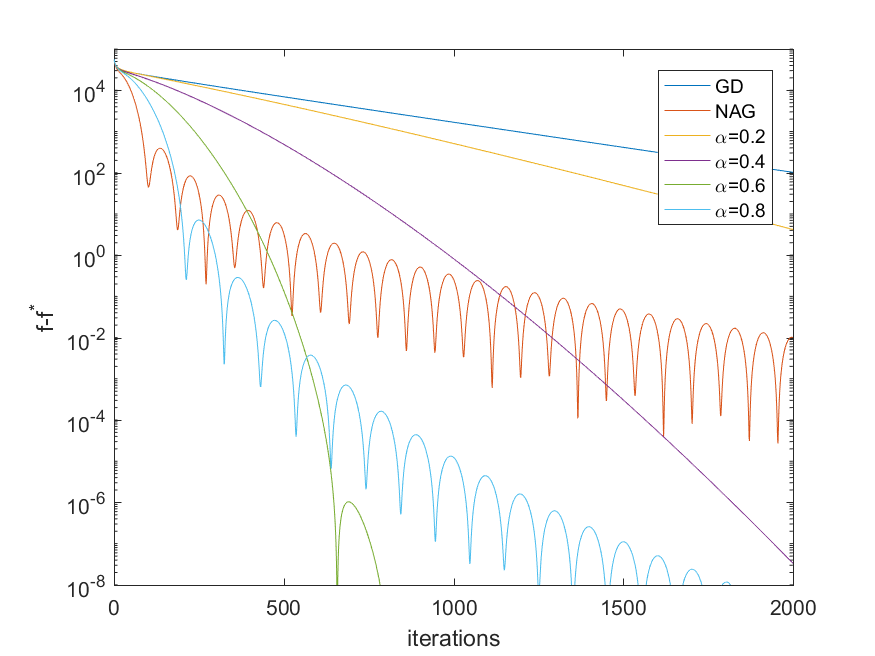}
  \end{minipage}
   }
  \subfigure[min \text{$\rho \log [\sum_{i=1}^m \exp((a_i^Tx-b_i)/\rho)]$}
  ]{\label{fig:log_a}
    \begin{minipage}[b]{0.45\textwidth}
    \includegraphics[width=\textwidth]{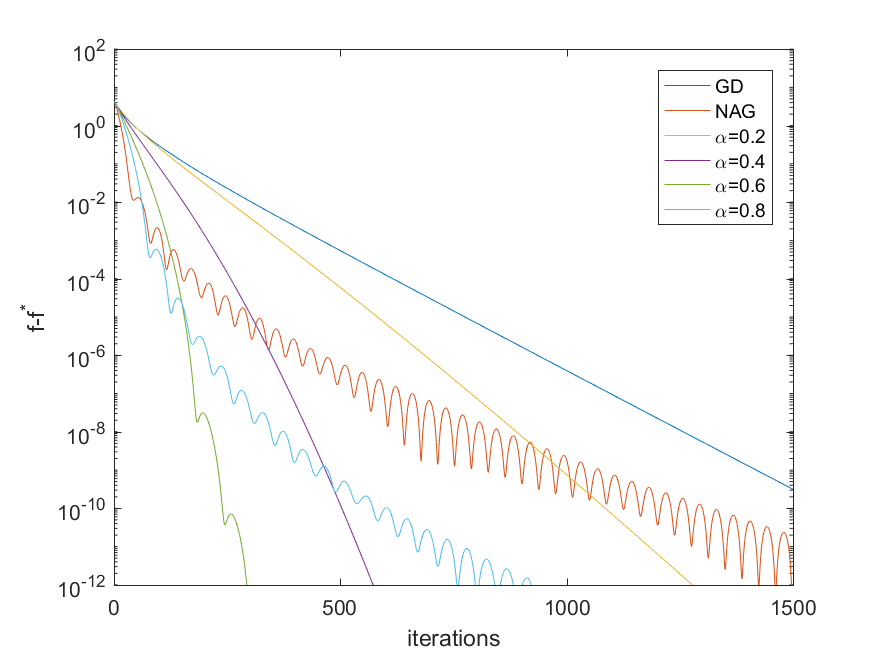}
  \end{minipage}
       }
\caption{Scheme (\ref{alg:SE_ar}) with fixed $r=3$ and varied $\alpha$}
\end{figure}

We take $GD$ (\ref{alg:GD}) and $NAG$ (\ref{alg:NAG}) as references. As shown in Figure \ref{fig:qudratic_a} and \ref{fig:log_a}, the symplectic Euler schemes (\ref{alg:SE_ar}) with different $\alpha$ have better numerical performances than GD and NAG when the number of iterations is large. Combining Theorem \ref{thm:convergence_ar} with Remark \ref{rmk:expression:T}, we know that although the convergence rate is faster as $\alpha$ gets smaller, the time $T$ when the convergence rate begins to be reached is relatively late. Therefore, for a trade-off between the convergence rate and the time $T$, we can see that (\ref{alg:SE_ar}) with $\alpha=0.6$ performs best among the four cases when $r=3$. NAG outperforms other algorithms in the early stage, and it is well-known that oscillations of NAG slow down its convergence rate and make it behave like GD. Based on the fact that the equation (\ref{eqn:E-L_ar}) is asymptotic to NAG's continuous analog (\ref{eqn:NAG}) as $\alpha$ approaches 1, the behavior of (\ref{alg:SE_ar}) is asymptotic to that of $NAG$ as $\alpha$ approaches 1. It follows that in the early stage, scheme (\ref{alg:SE_ar}) converges faster as $\alpha$ becomes larger. Besides, we can observe that oscillations become more frequent when $\alpha=0.8$ than $\alpha=0.6$, and there is no obvious oscillation when $\alpha=0.2$ and $\alpha=0.4$.\par
\begin{figure} 
\centering
\subfigure[min \text{$\frac{1}{2}x^TAx+b^Tx$}
]{ \label{fig:qudratic_r_a0.2}
\begin{minipage}[b]{0.45\textwidth}
    
    \includegraphics[width=\textwidth]{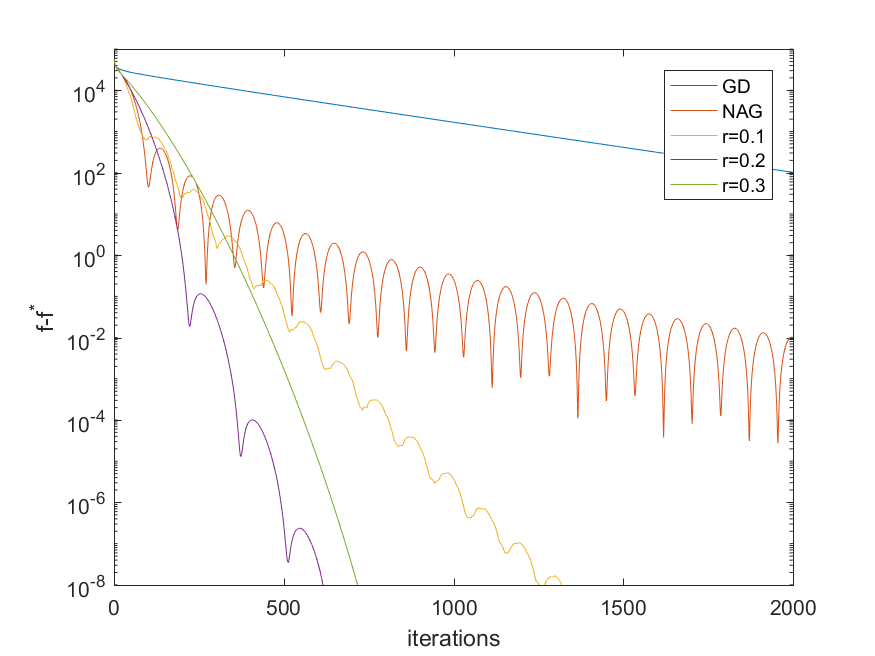}
  \end{minipage}
   }\vspace{-4mm}
   \subfigure[min \text{$\rho \log [\sum_{i=1}^m \exp((a_i^Tx-b_i)/\rho)]$}]{\label{fig:log_r_a0.2}
 
    \begin{minipage}[b]{0.45\textwidth}
    \includegraphics[width=\textwidth]{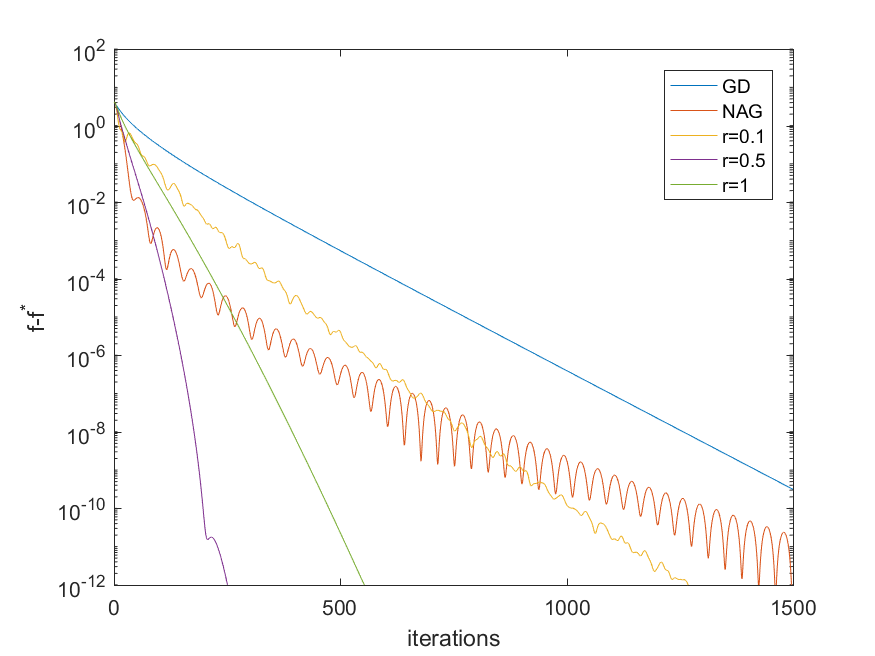}
  \end{minipage}
       } 
\subfigure[min \text{$\frac{1}{2}x^TAx+b^Tx$}]{ \label{fig:qudratic_r_a0.4}
\begin{minipage}[b]{0.45\textwidth}
    \includegraphics[width=\textwidth]{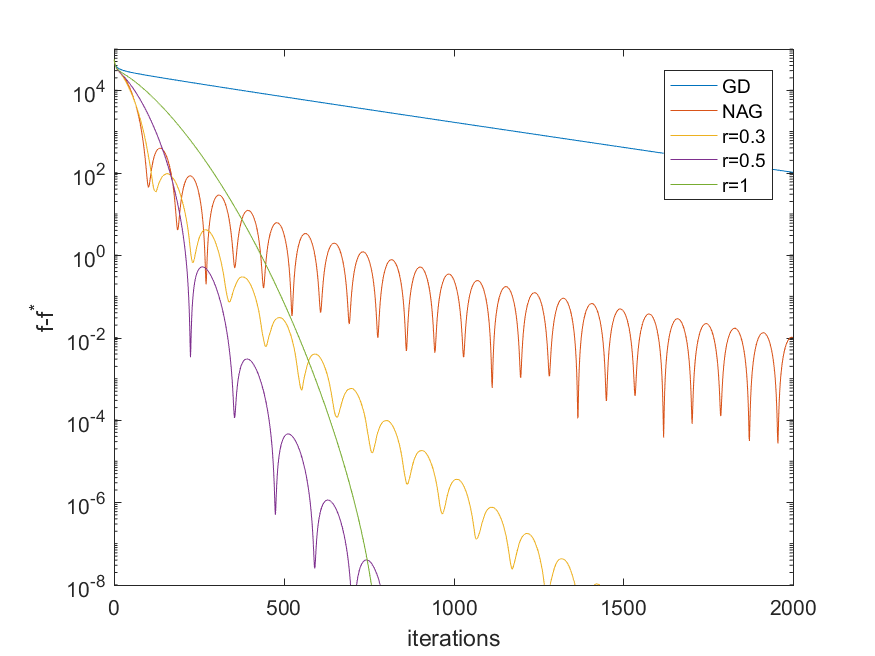}
  \end{minipage}
   } \vspace{-4mm}
\subfigure[min \text{$\rho \log [\sum_{i=1}^m \exp((a_i^Tx-b_i)/\rho)]$}]{\label{fig:log_r_a0.4}
    \begin{minipage}[b]{0.45\textwidth}
    \includegraphics[width=\textwidth]{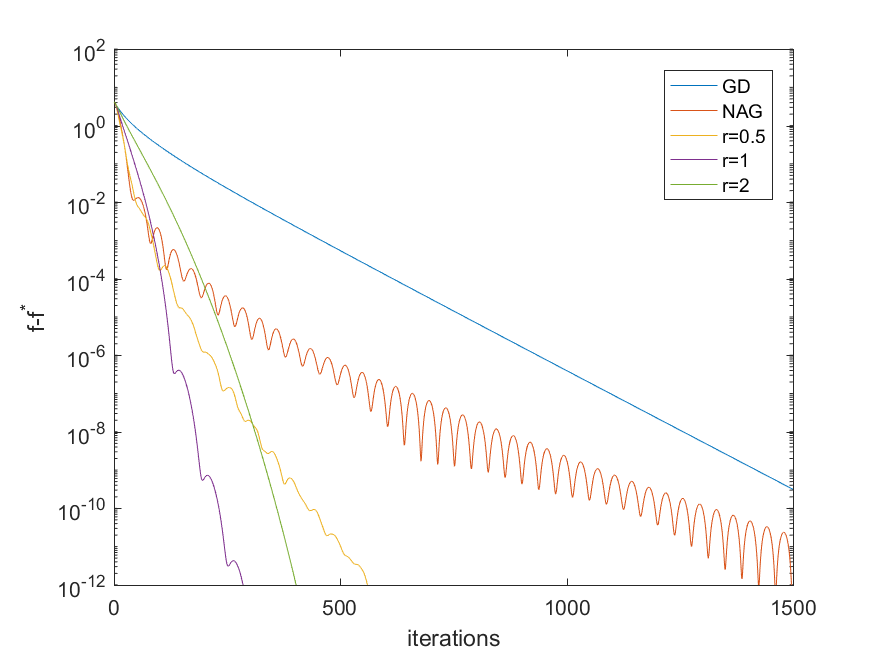}
  \end{minipage}
       } \subfigure[min \text{$\frac{1}{2}x^TAx+b^Tx$}]{ \label{fig:qudratic_r_a0.6}
\begin{minipage}[b]{0.45\textwidth}
    \includegraphics[width=\textwidth]{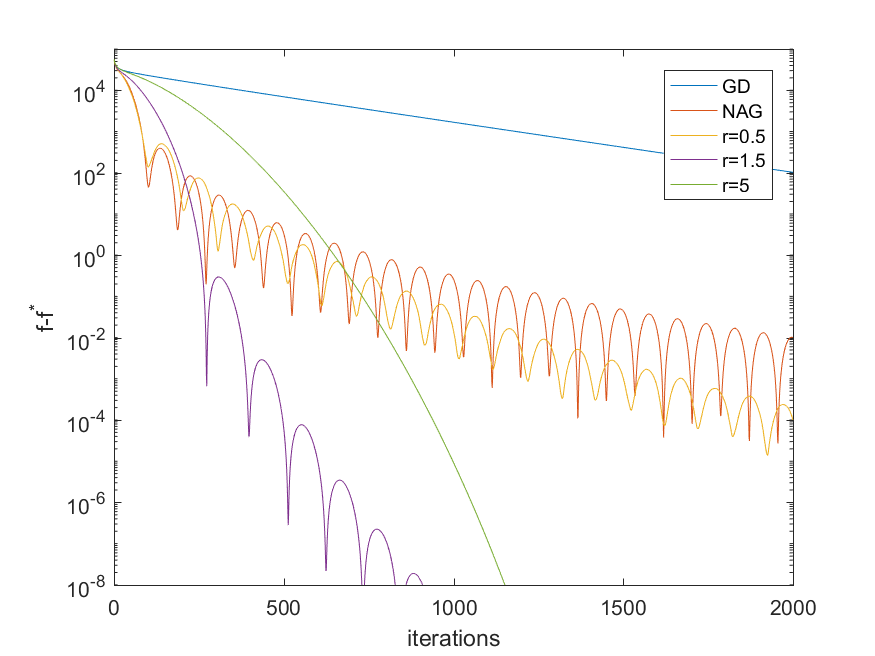}
  \end{minipage}
   }\vspace{-4mm}
\subfigure[min \text{$\rho \log [\sum_{i=1}^m \exp((a_i^Tx-b_i)/\rho)]$}]{\label{fig:log_r_a0.6}
    \begin{minipage}[b]{0.45\textwidth}
    \includegraphics[width=0.95\textwidth]{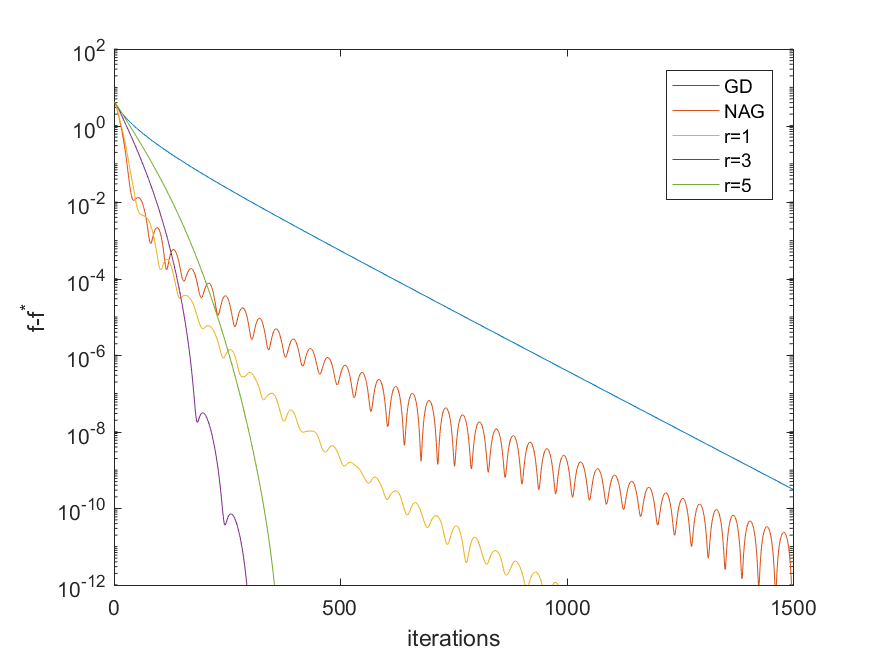}
  \end{minipage}
       } 
\subfigure[min \text{$\frac{1}{2}x^TAx+b^Tx$}]{ \label{fig:qudratic_r_a0.8}
\begin{minipage}[b]{0.45\textwidth}
    \includegraphics[width=\textwidth]{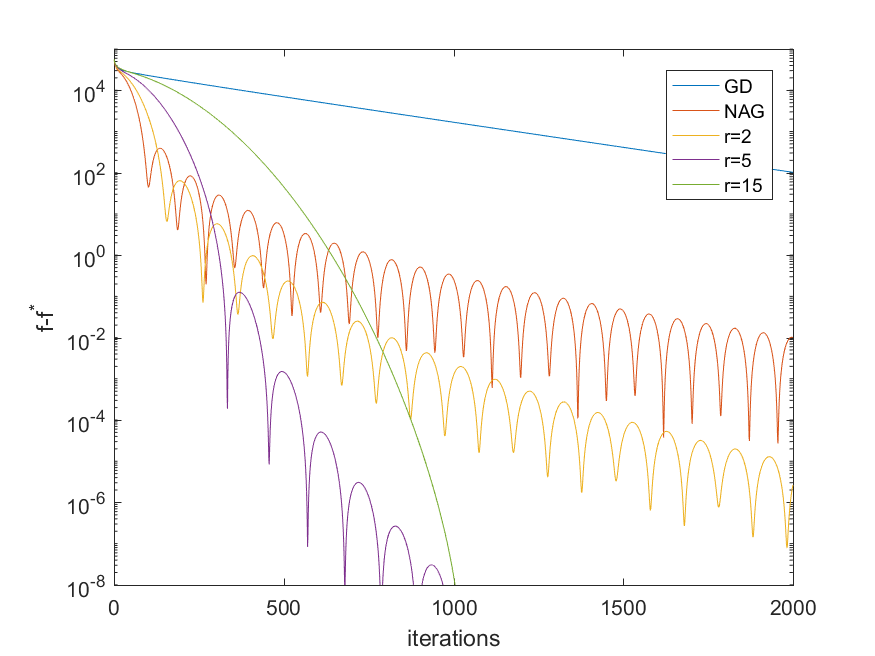}
  \end{minipage}
   } 
\subfigure[min \text{$\rho \log [\sum_{i=1}^m \exp((a_i^Tx-b_i)/\rho)]$}]{\label{fig:log_r_a0.8}
    \begin{minipage}[b]{0.45\textwidth}
    \includegraphics[width=\textwidth]{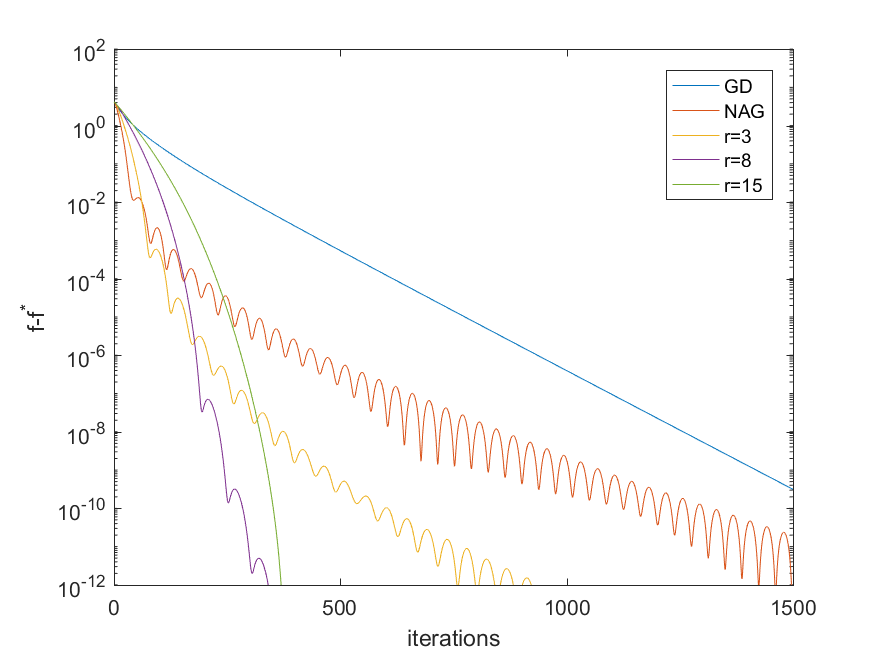}
  \end{minipage}
       }
\caption{Scheme (\ref{alg:SE_ar}) with fixed $\alpha=0.2,0.4,0.6,0.8$ and varied $r$}
\end{figure}
\begin{figure} [t] 
\centering
\subfigure[min \text{$\frac{1}{2}x^TAx+b^Tx$}
]{ \label{fig:qudratic_bestmatch}
\begin{minipage}[b]{0.45\textwidth}
    \includegraphics[width=\textwidth]{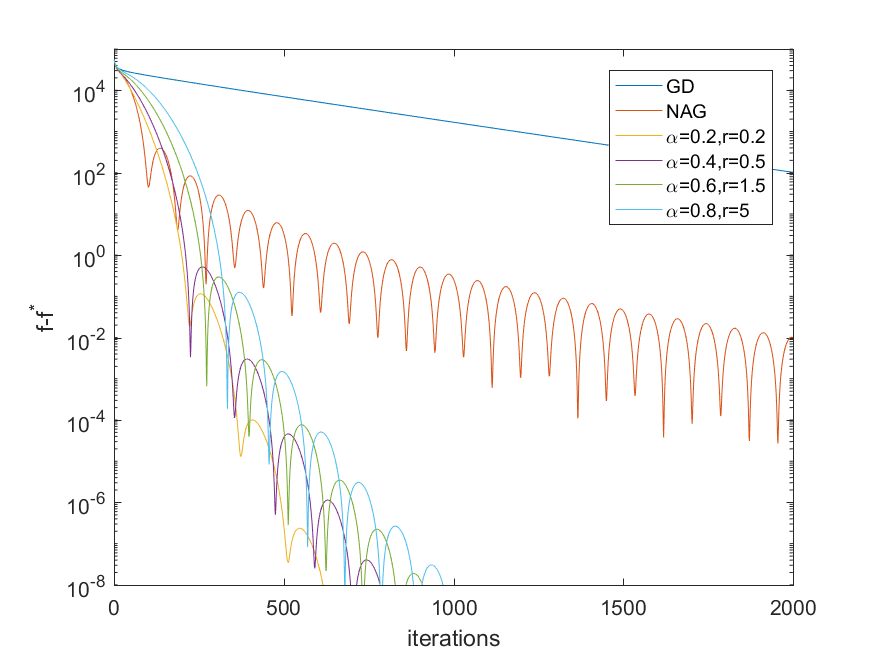}
  \end{minipage}
   }
  \subfigure[min \text{$\rho \log [\sum_{i=1}^m \exp((a_i^Tx-b_i)/\rho)]$}
  ]{\label{fig:log_bestmatch}
    \begin{minipage}[b]{0.45\textwidth}
    \includegraphics[width=\textwidth]{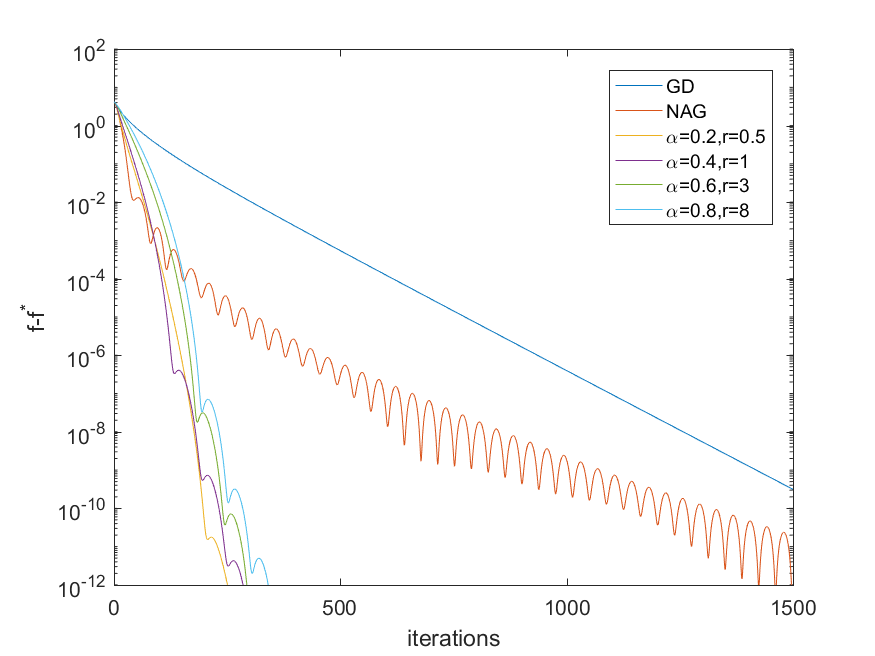}
  \end{minipage}
       }
\caption{Scheme (\ref{alg:SE_ar}) with best match ($\alpha,r$)}
\end{figure}

For fixed $\alpha$, we can see that the schemes (\ref{alg:SE_ar}) don't necessarily have better performance when $r$ gets larger. $r$ cannot be too large or too small to balance the convergence rate and the time $T$ from Theorem \ref{thm:convergence_ar} and Remark \ref{rmk:expression:T}. As depicted in Figure \ref{fig:qudratic_r_a0.2}, when $r=0.2$, it is better than when $r=0.1$ and $r=0.3$. However, the optimal $r$ increases as $\alpha$ increases. Let $r_{optimal}(\alpha)$ denotes the optimal $r$ for $\alpha$. For quadratic functions, $r_{optimal}(0.2)\approx 0.2,\ r_{optimal}(0.4)\approx 0.5,\ r_{optimal}(0.6)\approx 1.5$ and $r_{optimal}(0.8)\approx 5$ as illustrated in Figure \ref{fig:qudratic_r_a0.2}, \ref{fig:qudratic_r_a0.4}, \ref{fig:qudratic_r_a0.6} and \ref{fig:qudratic_r_a0.8}. Moreover, we can observe that when $r$ gets smaller, scheme (\ref{alg:SE_ar}) converges faster at the beginning and the oscillations becomes more frequent. In Figure \ref{fig:qudratic_r_a0.6}, the oscillation period of $r=0.5$ is similar as NAG but with smaller amplitude, while no oscillation appears if $r=5$. For different $\alpha$, we take the nearly optimal $r$. It follows that schemes (\ref{alg:SE_ar}) exhibit similar behavior as shown in Figure \ref{fig:qudratic_bestmatch} and \ref{fig:log_bestmatch}. 

\section{Conclusions and Outlook}\label{sect:5}
We know that the \textit{Bregman} dynamical system is a generalization of Nesterov's scheme in the continuous-time perspective. The \textit{Hamiltonian} of such system can be viewed as a \textit{separable} \textit{Hamiltonian} but with time-dependent kinetic energy and potential energy. We can construct integrators to approximate the solutions by constructing symplectic integrators of the extended autonomous \textit{Hamiltonian} system and projecting onto the original space. Since the extended \textit{Hamiltonian} is "separable" (we use the quotation since it is not separable for $t$), we apply the splitting methods, which are quite useful since they can preserve symplecticity like the exact flow and are usually easy to compute in practice. \par
An important characteristic of a numerical integrator is stability. Intuitively, the numerical orbit generated by an integrator should remain bounded if the exact orbit does. For dissipative system in this work, the exact solution of it is bounded. Therefore, we require the numerical solution to be bounded. In Section \ref{sect:2}, we get the linear stability region of scheme (\ref{alg:SE_xyuk}). Motivated by the fact that the \textit{Bregman} dynamical system can not be integrated via a single-loop scheme with a fixed step-size, we propose another kind of generalized Nesterov's scheme: ($\alpha,r$)-damped system. \par
In Section \ref{sect:3}, we construct Lyapunov functions and prove the convergence rate of the ($\alpha,r$)-damped system for strongly convex functions. These results provide us with perspective on the situation between the two Nesterov's schemes and offer insights into the choice of parameter $r$. In Lemma \ref{lemma:why strongly convex}, we illustrate why $f$ cannot just be convex. But it's technical. We hope to find other methods to show similar convergence rates for convex functions under some mild conditions, since we can observe in Section \ref{sect:4} that the scheme (\ref{alg:SE_ar}) has similar numerical performance for convex functions as strongly convex functions. \par
Besides, the relation between the two parameters $\alpha$ and $r$ is worth studying, since we can achieve better convergence performance if we properly choose $r$ for different $\alpha$. We expect to give quantitative analysis of the optimal $r$ with respect to $\alpha$, while we only give some qualitative results in this paper. Moreover, we can consider constructing some variants of the scheme (\ref{alg:SE_ar}) to escape saddle points in the non-convex setting as in \cite{escape}.

\bibliographystyle{alpha}
\bibliography{main}

\begin{thebibliography}{WWJ16}

\bibitem[Arn89]{arnol2013mathematical}
V.I. Arnold.
\newblock {\em Mathematical Methods of Classical Mechanics}.
\newblock Springer-Verlag, New York, second edition, 1989.

\bibitem[BCM12]{blanes2012splitting}
Sergio Blanes, Fernando Casas, and Ander Murua.
\newblock Splitting methods in the numerical integration of non-autonomous dynamical systems.
\newblock {\em Revista de la Real Academia de Ciencias Exactas, Fisicas y Naturales. Serie A. Matematicas}, 106:49--66, 2012.

\bibitem[BJW18]{betancourt2018symplectic}
Michael Betancourt, Michael~I Jordan, and Ashia~C Wilson.
\newblock On symplectic optimization.
\newblock {\em arXiv preprint arXiv:1802.03653}, 2018.

\bibitem[CSY22]{chen2022gradient}
Shuo Chen, Bin Shi, and Ya-xiang Yuan.
\newblock Gradient norm minimization of {Nesterov} acceleration: $ o (1/k^{3}) $.
\newblock {\em arXiv preprint arXiv:2209.08862}, 2022.

\bibitem[DJ21]{diakonikolas2021generalized}
Jelena Diakonikolas and Michael~I Jordan.
\newblock Generalized momentum-based methods: A {Hamiltonian} perspective.
\newblock {\em SIAM Journal on Optimization}, 31(1):915--944, 2021.

\bibitem[FQ10]{feng2010symplectic}
Kang Feng and Mengzhao Qin.
\newblock {\em Symplectic Geometric Algorithms for Hamiltonian Systems}.
\newblock Zhejiang Science and Technology Publishing House, Hangzhou and Springer-Verlag Berlin Heidelberg, 2010.

\bibitem[HLW06]{hairer2006geometric}
Ernst Hairer, Christian Lubich, and Gerhard Wanner.
\newblock {\em Geometric Numerical Integration: Structure-Preserving Algorithms for Ordinary Differential Equations}.
\newblock Springer-Verlag Berlin Heidelberg, second edition, 2006.

\bibitem[JNJ18]{escape}
Chi Jin, Praneeth Netrapalli, and Michael~I Jordan.
\newblock Accelerated gradient descent escapes saddle points faster than gradient descent.
\newblock In {\em Conference On Learning Theory}, pages 1042--1085. PMLR, 2018.

\bibitem[KY23]{kim2023unifying}
Jungbin Kim and Insoon Yang.
\newblock Unifying {Nesterov’s} accelerated gradient methods for convex and strongly convex objective functions.
\newblock In {\em International Conference on Machine Learning}, pages 16897--16954. PMLR, 2023.

\bibitem[LC22]{luo2022differential}
Hao Luo and Long Chen.
\newblock From differential equation solvers to accelerated first-order methods for convex optimization.
\newblock {\em Mathematical Programming}, 195(1):735--781, 2022.

\bibitem[LSY24]{li2024linear}
Bowen Li, Bin Shi, and Ya-xiang Yuan.
\newblock Linear convergence of forward-backward accelerated algorithms without knowledge of the modulus of strong convexity.
\newblock {\em SIAM Journal on Optimization}, 34(2):2150--2168, 2024.

\bibitem[Nes83]{nesterov1983method}
Y~Nesterov.
\newblock A method for solving a convex programming problem with convergence rate ${O} (1/k^2)$.
\newblock In {\em Soviet Mathematics. Doklady}, volume~27, pages 367--372, 1983.

\bibitem[Nes04]{nesterov2004introductory}
Yurii Nesterov.
\newblock {\em Introductory Lectures on Convex Optimization: A Basic Course}, volume~87.
\newblock Springer Science \& Business Media, 2004.

\bibitem[Per01]{perko2013differential}
Lawrence Perko.
\newblock {\em Differential Equations and Dynamical Systems}.
\newblock Springer-Verlag, New York, third edition, 2001.

\bibitem[SBC16]{NAGODE}
Weijie Su, Stephen Boyd, and Emmanuel~J Candes.
\newblock A differential equation for modeling {Nesterov's} accelerated gradient method: Theory and insights.
\newblock {\em Journal of Machine Learning Research}, 17(153):1--43, 2016.

\bibitem[SDJS22]{shi2022understanding}
Bin Shi, Simon~S Du, Michael~I Jordan, and Weijie~J Su.
\newblock Understanding the acceleration phenomenon via high-resolution differential equations.
\newblock {\em Mathematical Programming}, pages 1--70, 2022.

\bibitem[WRJ21]{wilson2021lyapunov}
Ashia~C Wilson, Ben Recht, and Michael~I Jordan.
\newblock A {Lyapunov} analysis of accelerated methods in optimization.
\newblock {\em Journal of Machine Learning Research}, 22(113):1--34, 2021.

\bibitem[WWJ16]{1603}
Andre Wibisono, Ashia~C Wilson, and Michael~I Jordan.
\newblock A variational perspective on accelerated methods in optimization.
\newblock {\em proceedings of the National Academy of Sciences}, 113(47):E7351--E7358, 2016.

\end{thebibliography}

\end{document}